\documentclass[twoside,11pt]{article}
\usepackage[preprint]{jmlr2e}
\usepackage{amsmath,mathtools,booktabs,array,enumitem,xcolor}
\usepackage{tikz}
\usetikzlibrary{trees}
\usepackage[ruled,vlined,linesnumbered]{algorithm2e}

\allowdisplaybreaks[2]
\numberwithin{equation}{section}
\newtheorem{assumption}[theorem]{Assumption}

\newcommand{\norm}[1]{\lVert#1\rVert}
\newcommand{\abs}[1]{\lvert#1\rvert}
\newcommand{\ip}[2]{\langle#1,#2\rangle}
\newcommand{\ee}{\mathrm e}
\newcommand{\dd}{\,\mathrm d}
\newcommand{\Acal}{\boldsymbol a}
\newcommand{\cT}{\mathcal T}
\newcommand{\cU}{\mathcal U}

\SetKwInput{KwIn}{Input}
\SetKw{Return}{return}
\SetKw{Break}{break}
\SetKw{Continue}{continue}
\DontPrintSemicolon

\hypersetup{
	pdftitle={Adaptive Chebyshev Gradient Descent},
	pdfauthor={Yin Liu, Ziyang Zeng, and Junyu Zhang},
	bookmarksnumbered=true
}
\ShortHeadings{Adaptive Chebyshev Gradient Descent}{Liu, Zeng, and Zhang}
\firstpageno{1}

\begin{document}
	
	\title{Near Optimal Pure Stepsize Adaptation Acceleration for \\Gradient Descent on Nonconvex Optimization }
	\author{\name Yin Liu \email yinliu@pku.edu.cn \\
		\addr Beijing International Center for Mathematical Research\\ Peking University
		\AND
		\name Ziyang Zeng \email ziyangzeng@u.nus.edu \\
		\addr Department of Industrial Systems Engineering and Management\\ National University of Singapore
		\AND
		\name Junyu Zhang \email junyuz@nus.edu.sg \\
		\addr Department of Industrial Systems Engineering and Management\\ National University of Singapore}
	
	\maketitle
	
	\begin{abstract}%
		Existing accelerated first-order methods for nonconvex optimization typically rely on momentum or extrapolation, requiring auxiliary vectors beyond the current iterate and gradient. On the other hand, memory-free stepsize scheduling alone can accelerate certain convex problems, while no guarantee exists for nonconvex problems. This paper closes this gap by proposing Nonconvex Adaptively Clipped Chebyshev (NACC) gradient descent method, which achieves near optimal acceleration while keeping the vanilla gradient descent update. NACC combines an adaptive selection of reversed fractal Chebyshev stepsizes with an epoch-wise clipping mechanism to control the aggressiveness of stepsizes. For nonconvex problems with Lipschitz continuous gradients and Hessians, we establish a near-optimal $\widetilde O(\epsilon^{-7/4})$  complexity. If the problem is convex in addition, then, without modification, NACC automatically achieves an $\widetilde O(\epsilon^{-3/4})$ accelerated complexity. Importantly, NACC guarantees the standard quadratic sufficient descent at every iteration, even when the stepsizes are substantially larger than the conventional choice.   These monotonicity properties further facilitate practical implementations that adaptively adjust the estimated gradient and Hessian Lipschitz constants, demonstrating significant advantages in numerical tests against the existing nonconvex acceleration benchmarks. 
	\end{abstract}
	
	
\section{Introduction}
Accelerating first-order methods for minimizing a nonconvex function $\min_x f(x)$ has been an important topic for years. Early results attempt to extend Nesterov's accelerated gradient descent method to nonconvex problems, see e.g., \citep{ghadimi2016accelerated,ghadimi2019generalized}. However, for finding an $\epsilon$-stationary point, a point $\bar x$ where $\|\nabla f(\bar x)\|\leq \epsilon$, under the general $L$-smooth setting, such methods are not able to improve the $O(\epsilon^{-2})$ first-order oracle complexity of the pure gradient descent (GD) method. In fact, the lower bounds established by \citet{carmon2020lower,carmon2021lower} show that this complexity order is unimprovable under the standard $L$-smooth assumption alone. Faster results are possible only when exploiting additional structure. For example, given both gradient and Hessian Lipschitz continuity, which is also the main problem setting that this paper focuses on, \citet{carmon2021lower} established an $\Omega(\epsilon^{-12/7})$   oracle complexity lower bound for first-order methods, and is recently tightened to $\Omega(\epsilon^{-7/4})$ by \citet{zhou2026sharp}. 

Indeed, under the gradient and Hessian Lipschitz continuity, \citet{carmon2017convex} developed a convex until proven guilty (CUPG) framework and achieved $\widetilde O(\epsilon^{-7/4})$ complexity using only first-order oracles. CUPG applies a monitored accelerated gradient descent to proximal point type subproblems, achieving either a convex-like accelerated progress or a negative curvature direction that leads to a sufficient descent. Nevertheless, CUPG has a sophisticated structure of multiple nested subroutines, making it hard to implement and tune in practice. 
Subsequent works by \citet{li2023restarted} and \citet{marumo2024parameterfree} focused on improving the simplicity and practical implementability of accelerated first-order methods under this setting, achieving the same $O(\epsilon^{-7/4})$ complexity without negative curvature search or proximal subproblems. 

Despite these advances, existing accelerated methods still rely on momentum or extrapolation, which requires maintaining auxiliary vectors beyond the current iterate $x_k$ and gradient $\nabla f(x_k)$. As the scale of modern optimization problems grows, storing and updating these vectors incurs additional memory and per-iteration computational costs. This motivates us to consider whether the same near-optimal complexity can be achieved with a simpler iteration structure. In addition, beyond the practical memory and computation efficiencies, this also raises a more fundamental question about the mechanism of nonconvex acceleration, which is also the main focus of this paper: 

\begin{quote}
    \emph{For a nonconvex optimization problem with Lipschitz continuous gradient and Hessian, can we achieve the (near) optimal accelerated rate only via stepsize adaptation, while maintaining the simple GD update form? }
\end{quote}

A positive clue comes from the existing works on convex optimization. For general $L$-smooth convex problems,  several recent studies \citep{altschuler2025acceleration,grimmer2025accelerated,zhang2026accelerated} achieved a (non-optimal) $O(\epsilon^{-0.786})$ accelerated complexity for GD via carefully designed stepsize schedules. For convex quadratic functions,  a related line of work \citep{young1953richardson,golub1961chebyshev,lebedev1971ordering,agarwal2021acceleration} developed a Chebyshev stepsize schedule, which, interestingly, achieved the optimal $O(\sqrt{\kappa}\log(1/\epsilon))$ complexity for the strongly convex quadratic case, with $\kappa$ denoting the condition number of the problem. However, neither of these methods converge for nonconvex problems. For Chebyshev schedule, despite its optimality for strongly convex quadratics, it may not even converge for non-quadratic strongly convex problems, see \citep{agarwal2021acceleration}. To the best of our knowledge, there is no existing results that directly accelerate GD for nonconvex optimization by pure stepsize scheduling.  

\paragraph{Our contributions.} We answer our question affirmatively by proposing the Nonconvex Adaptively Clipped Chebyshev (NACC) gradient descent method. It takes the pure GD update form $x_{k+1}=x_k-\eta_k\nabla f(x_k)$ while only adjusting the stepsize $\eta_k$ through a two-level adaptive mechanism based on the reversed fractal Chebyshev schedule. At the beginning of each epoch, NACC adaptively selects a block of Chebyshev stepsizes by the current gradient size and smoothness constants, determining the potential aggressiveness of this stepsizes in the epoch. Within the epoch, an adaptive clipping mechanism is imposed to further control the actual executed stepsizes according to the obtained progress. Based on this method, we obtain the following results.

\begin{itemize}[leftmargin=0.5cm]
    \item \textbf{Stepsize-only nonconvex acceleration.} We design the NACC algorithm, the first method to accelerate GD in nonconvex optimization through pure stepsize adaptation.
    \item \textbf{Near-optimal nonconvex complexity.} Assuming the Lipschitz continuity of gradient and Hessian, NACC returns an $\epsilon$-stationary point within $\widetilde O(\epsilon^{-7/4})$ gradient evaluations, matching the optimal lower bound up to a square root of logarithmic factor. 
    \item \textbf{Automatic convex acceleration.} Without any modification, NACC automatically achieves an accelerated complexity of $\widetilde O(\epsilon^{-3/4})$ for convex objectives, in terms of finding an $\epsilon$-optimal point $\bar x$ such that $f(x)-f^*\le \epsilon$.
    \item \textbf{Monotonicity.} Unlike many accelerated gradient methods, NACC is a monotonic acceleration method.  Regardless of convexity, it always satisfies a quadratic sufficient descent property:  $f(x_{k+1})\leq f(x_k)-c\eta_k\|\nabla f(x_k)\|^2$. For convex problems, we further have $\|x_{x+1}-x^*\|\leq \|x_k-x^*\|$.  
\end{itemize}

In particular, it is worth mentioning that the monotonicity properties provides the potential of our algorithm for stronger convergence guarantees under additional function structures, such as the sequence convergence of tame objective functions \citep{attouch2013convergence,bolte2014proximal}. In addition, the sufficient descent property facilitates a more flexible descent-based adaptive implementation, where the estimates of the gradient and Hessian Lipschitz constants can be adjusted simultaneously according to the observed function decrease, based on which a practically implemented NACC with adaptively adjusted Lipschitz constants consistently demonstrate the existing nonconvex acceleration benchmarks \citep{carmon2017convex,li2023restarted,marumo2024parameterfree}.

\paragraph{Comparison with existing Chebyshev stepsize results.} Finally, we would like to end the introduction with a comparison with the most related works, the existing Chebyshev stepsize schedule results \citep{young1953richardson,golub1961chebyshev,lebedev1971ordering,agarwal2021acceleration}. The fundamental difference between our results and existing Chebyshev stepsize methods lies in the transition from convex quadratics to general nonconvex non-quadratic objectives, which introduces two major challenges. First, the negative curvature of nonconvex problems destroys the stability of Chebyshev schedule, which is a key tool for existing works. Second, the Hessian variation in non-quadratic objectives leads to non-commuting matrix products, another key tool for existing works.  To overcome these difficulties, our NACC method replaces the predetermined Chebyshev schedules with a two-level adaptive Chebyshev schedule, and eventually retains the near optimal complexity of Hessian-Lipschitz nonconvex optimization, as opposed to the existing results that may not even converge for strongly convex problems \cite{agarwal2021acceleration}.

\section{Preliminaries on Reversed Fractal Chebyshev Schedule} 
\label{sec:prelim}

An important component of our adaptive stepsize rule is the  Chebyshev schedule that constructs the stepsizes from the roots of Chebyshev polynomials \citep{young1953richardson,golub1961chebyshev,lebedev1971ordering,agarwal2021acceleration}. These results focus on strongly convex quadratics $f(x)=\frac{1}{2}x^\top H x-b^\top x$ with $\mu I\preceq H \preceq LI.$ As taking $K$ GD updates for this problem with stepsizes $\eta_0,\ldots,\eta_{K-1}$ leads to a linear dynamics of
$$
x_K-x^\star = \prod_{k=0}^{K-1}(I-\eta_k H) (x_0-x^\star),
$$
the fastest stepsize schedule then corresponds to selecting $\eta_t, t\leq K-1$ that minimizes  $\|\Pi_{k=0}^{K-1}(I-\eta_k H)\|$. For such a product of commutative matrices, it reduces to minimizing $\max_{\lambda\in[\mu,L]}|\Pi_{k=0}^{K-1}(1-\eta_k\lambda)|$, whose solution is exactly given by the roots of some affinely scaled Chebyshev polynomial \citep{golub1961chebyshev}. Later on,  \citep{lebedev1971ordering} designed a specific ordering of the roots, called fractal Chebyshev schedule, and showed that it can significantly improve the method's stability. This ordering has recently been revisited by \citep{agarwal2021acceleration}, in which a reversed variant of this ordering is analyzed. Due to the three-term recurrence of the Chebyshev polynomial, the last point $x_K$ provably coincides with a momentum-based update sequence, though the intermediate iterates are generally different, see \citep{daspremont2021acceleration}
for details of this interesting connection with momentum methods.

Back to our nonconvex and non-quadratic setting, in this paper, we mainly adopt the reversed fractal ordering, which can be defined recursively as
	\begin{equation}
		\label{eq:A-def}
		\Acal_{1,\phi} = \bigl[p_{1,\phi}\bigr], \quad\Acal_{2m,\phi}=\bigl[\Acal_{m,\pi-\phi/2},\ \Acal_{m,\phi/2}\bigr], \quad
		p_{m,\phi} := \csc^2\!\left(\frac{\phi}{2m}\right),
	\end{equation}
	where $m=2^k$ for $k=0,1,2\cdots$,  $\phi\in(0,\pi)$, and $\csc(\theta):=1/\sin(\theta)$ stands for the cosecant function.  It is worth mentioning that our reversed fractal scheduling \eqref{eq:A-def} is different from   \cite{agarwal2021acceleration}, as we do not rescale the polynomial to some spectral interval $[\mu,L]$ and we allow a general parameter $\phi$ to facilitate the analysis of nonconvex problems.
	\begin{figure}[ht]
		\centering
		\begin{tikzpicture}[
			level distance=1.2 cm,
			level 1/.style={sibling distance=7.2cm},
			level 2/.style={sibling distance=3.4cm},
			level 3/.style={sibling distance=1.7cm},
			every node/.style={align=center},
			phase/.style={font=\small},
			rootval/.style={font=\small}
			]
			\node[phase] {$\Acal_{4,\phi}$}
			child {
				node[phase] {$\Acal_{2,\pi-\phi/2}$}
				child {
					node[phase] {$\Acal_{1,\pi/2+\phi/4}$}
					child {
						node[rootval] {$\displaystyle \csc^2\!\left(\frac{\pi}{4}+\frac{\phi}{8}\right)$}
					}
				}
				child {
					node[phase] {$\Acal_{1,\pi/2-\phi/4}$}
					child {
						node[rootval] {$\displaystyle \csc^2\!\left(\frac{\pi}{4}-\frac{\phi}{8}\right)$}
					}
				}
			}
			child {
				node[phase] {$\Acal_{2,\phi/2}$}
				child {
					node[phase] {$\Acal_{1,\pi-\phi/4}$}
					child {
						node[rootval] {$\displaystyle \csc^2\!\left(\frac{\pi}{2}-\frac{\phi}{8}\right)$}
					}
				}
				child {
					node[phase] {$\Acal_{1,\phi/4}$}
					child {
						node[rootval] {$\displaystyle \csc^2\!\left(\frac{\phi}{8}\right)$}
					}
				}
			};
		\end{tikzpicture}
		\caption{Recursive construction of $\Acal_{4,\phi}$, with tree height $h_m = 1+\log_2m$ for $ m=4$.  }
		\label{fig:tree-m4}
	\end{figure}
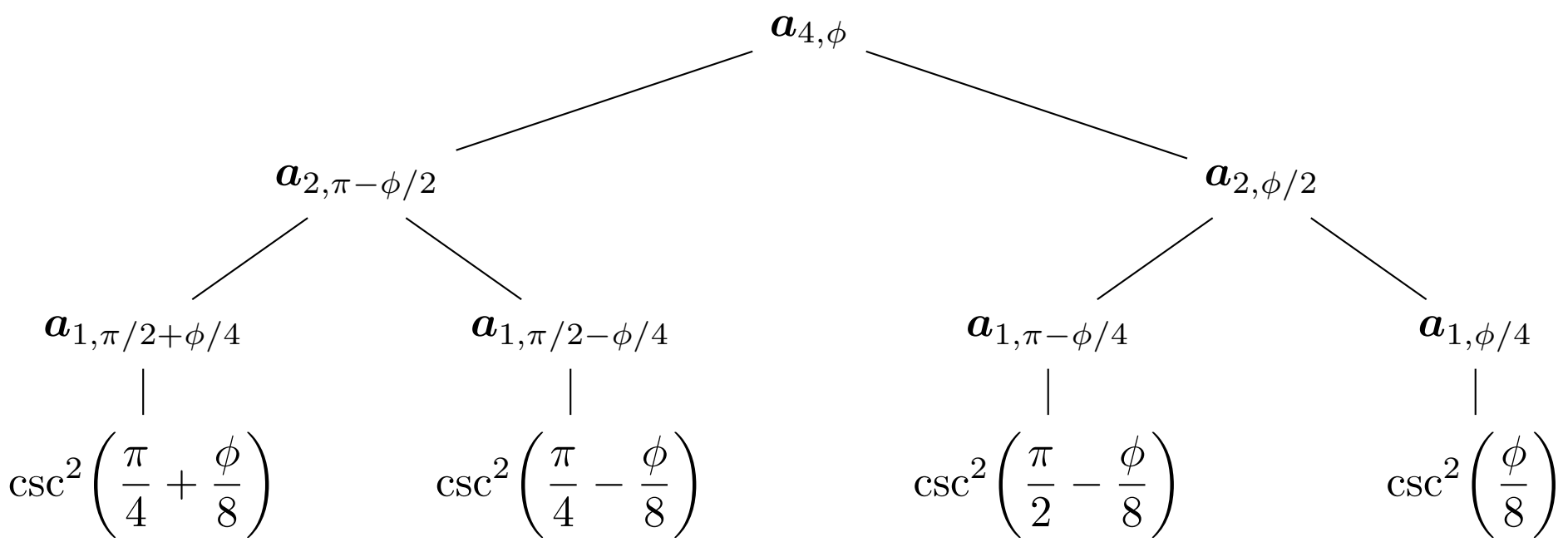
	
	For the ease of understanding, we provide a detailed example for the construction of $\Acal_{4,\phi}$ in Figure \ref{fig:tree-m4}. As $\Acal_{m,\phi}$ can be represented as a binary tree, we denote its height as $h_m:=1+\log_2m$.  With a standard choice $\phi = \pi/2$,  Figure \ref{fig:tree-m4}  gives 
	$$\Acal_{4,\pi/2}  = \Big[\csc^2\!\big(\frac{5\pi}{16}\big), \  \csc^2\!\big(\frac{3\pi}{16}\big), \ \csc^2\!\big(\frac{7\pi}{16}\big), \   \csc^2\!\big(\frac{\pi}{16}\big)\Big].$$
	Longer sequences can be obtained recursively in a similar way. The sequence \eqref{eq:A-def} has rich properties. In detail, recall the Chebyshev polynomial of the first kind $\cT$. For a finite sequence $Q=[a_1,\ldots,a_{r}]$, define its associated polynomial by
	\begin{equation} 
		P_Q(z):=\prod_{i=1}^{r}(1-a_i z), \nonumber
	\end{equation}
	and default $P_{\varnothing}(z):=1$ for empty sequence. 
    With a direct rescaling of the skewed Chebyshev polynomial representation in \citet[Appendix C.1]{agarwal2021acceleration}, we obtain
	\begin{equation}\label{eq:PQ-cheb}
		P_{\Acal_{m,\phi}}(z)
		=\frac{\cT_m(1-2z)-\cos\phi}{1-\cos\phi},
	\end{equation}
    where $\cT$ denotes the Chebyshev polynomials of the first kind that is characterized by $\cT_n(\cos t)=\cos(nt)$. 
	It is also straightforward to verify by induction that the sequence $\Acal_{m,\phi}$ is a reordering of the set $\big\{\csc^2\!\big(\frac{2\pi j+\phi}{2m}\big)
	:j=0,\ldots,m-1
	\big\},$
	which consists of the (inversed) roots of the equation $\cT_m(1-2z)=\cos \phi$. Among this set, $p_{m,\phi} = \csc^2(\phi/2m)$ is the largest ``peak'' element, this is why we choose ``$p$'' in the notation.  Then, though not found in existing works to our best knowledge, it is not hard to obtain for $\phi\in[\pi/2,\pi)$ that 
	\begin{equation}\label{eq:general-mass}
		\sum_{a\in\Acal_{m,\phi}}a
		=\frac{2m^2}{1-\cos\phi}\in[m^2,2m^2] 
	\end{equation}
	by differentiating \eqref{eq:PQ-cheb} at $z=0$ and using the fact $\cT_m'(1)=m^2$. For the specific choice of $\phi = \pi/2$ in later discussion, we further have $m^2\leq p_{m,\pi/2}\leq2m^2$. 
	Finally, slightly generalizing the proof of \citet[Proposition 13]{agarwal2021acceleration} to general $\phi$ gives the stability lemma.
	\begin{lemma}
		\label{lemma:prefix-stability}
		Let $Q = \Acal_{m,\phi} = (a_1,\cdots,a_m)$, with $\phi\in[\pi/2,\pi)$ and $m=2^k$ for some integer $k\geq0$. Then for any $Q_{1:r} := (a_1,\cdots,a_r), r\leq m$, we have $|P_{Q_{1:r}}(z)|\leq 1, \forall z\in[0,1].$
	\end{lemma} 
	The reason why such a property is called \emph{stability} can be illustrated by recalling the convex quadratic setting discussed in the introduction. When we select the stepsize sequence  as $[\eta_1,\ldots,\eta_m]=\frac{1}{L}\cdot\Acal_{m,\phi}$, for some $L\geq\|H\|$ and $\phi\in[\pi/2,\pi)$, one can always guarantee 
	$$\|x_k-x^*\| = \|P_{Q_{1:k}}(H/L)(x_0-x^*)\|\leq \|x_0-x^*\|, \quad \forall k\leq m.$$
	Such a property does not always hold if this sequence of stepsizes is arbitrarily permuted, and hence we say the reversed fractal Chebyshev stepsize schedule is stable. 
    
    However, we should also note that for indefinite $H\not\succeq0$ that corresponds to a nonconvex quadratic problem, the eigenvalues of $H/L$ lies in $[-1,1]$ rather than $[0,1]$, and hence $\|P_{Q_{1:k}}(H/L)\|$ is no longer controlled by 1, and the stability will fail.

	In addition, as our objective function is neither convex nor quadratic, and the fixed Hessian matrix $H$ in the above naive explanation will change to a sequence of noncommutative matrices that vary along the iterations. Hence, we cannot adopt the direct Chebyshev stepsize scheme of \citet{agarwal2021acceleration}, and we cannot rely on the simple polynomial-based stability analysis. How to handle the negative curvatures from nonconvexity and the drifting Hessian issue from the non-quadratic objective therefore motivates the adaptive algorithm discussed in the next section.
	
	\section{Adaptively Clipped Chebyshev GD for Nonconvex Problems}
	\label{sec:nonconvex}
	In this section, we begin with the standard assumptions of our main problem and the corresponding notations, and then we introduce our Chebyshev GD algorithm with adaptive clipping and present its main theoretical results, with part of the heavy analysis relegated to the appendix.

	\subsection{Basic Problem Setting and the Proposed Algorithm} Throughout the paper, we will consider the unconstrained minimization of $f$ that satisfies the following second-order smooth and lower bounded assumption. 
	
	\begin{assumption}
		\label{assumption:f}
		$f:\mathbb{R}^d\to\mathbb{R}$ is twice continuously differentiable. There exist $L,M>0$ s.t.
		$$\|\nabla f(x)-\nabla f(y)\|
		\leq L\|x-y\| \quad \mbox{and} \quad \|\nabla^2 f(x) - \nabla^2 f(y)\| \leq M\|x-y\|$$  
		hold for all $x,y \in \mathbb{R}^d$. In addition, we assume $f^*:=\inf_x f(x)>-\infty$.
	\end{assumption}
	
	Under this assumption, the following standard descent lemma \citep{nesterov2006cubic} holds that
	\begin{equation}
		\label{eq:cubic} 
		f(x+\xi) \leq f(x) + \xi^\top\nabla f(x)  + \frac{1}{2}\xi^\top\nabla^2 f(x)\xi + \frac{M}{6}\|\xi\|^3. 
	\end{equation} 
	In particular, for any symmetric matrix $H$, denote $[H]_+\succeq0$ the projection of $H$ to the positive semidefinite cone. Then a handy inequality is $\xi^\top H\xi\leq \xi^\top[H]_+\xi$, for $\forall \xi$. 
	
	\subsection{Adaptively Clipped Chebyshev (NACC) Gradient Descent}
	
	Assuming prior knowledge of the Lipschitz constants $L$ and $M$, we present the following Nonconvex Adaptively Clipped Chebyshev (NACC) gradient descent method in Algorithm \ref{alg:nonconvex}. In this algorithm, the parameter $\phi$ of \eqref{eq:A-def} is defaulted at $\phi = \pi/2$ and the fixed constant parameter $\theta = 19/20<1$ as well as the $1/3$ factor in \eqref{eq:nonconvex-cert} are determined specifically to simplify the analysis under the current choice of $\phi$. In practice, they are tunable parameters. 
	
	\begin{algorithm}
		\caption{Adaptively Clipped Chebyshev (NACC) Gradient Descent}\label{alg:nonconvex}
		\small
		\textbf{Initialize:}\, Initial point $\widehat{x}_0$, Lipschitz constants $L,M>0$, tolerance $\epsilon>0$, and $t=0$, $\theta = \frac{19}{20}$\; \vspace{0.1cm}
		\While{$\|\nabla f(\widehat{x}_t)\|>\epsilon$}{ 
			Set $x_{t,0}= \widehat{x}_{t}$,  $G_t=\|\nabla f(\widehat{x}_t)\|$, and choose the largest $m\in\{4,8,16,\ldots\}$ that satisfies \eqref{eq:nonconvex-cert}.\;
			\If{no such $m$ exists}{
				Set $\widehat{x}_{t+1} = \widehat{x}_{t} - \frac{3}{2L}\nabla f(\widehat{x}_{t})$ and update $t\gets t+1$.\;
				\Continue.\;
			}
			Denote $\Acal_{m,\pi/2}=[a_{t,0},\ldots,a_{t,m-1}]$, and set 
			$i\gets0$, $E_{t,i}\gets0$, $E_{\max}^t\gets \frac{1}{2}{\theta^2G_t^2\cdot p_{m,\pi/2}}. $\;
			\While{$i<m$}{
				\textbf{if} $\norm{\nabla f(x_{t,i})}\leq\epsilon$ \textbf{then return} $\widehat{x}_{t+1} = x_{t,i}$.\; 
				
				\textbf{if} $\norm{\nabla f(x_{t,i})}\leq\theta G_t$ \emph{or} $E_{t,i}=E_{\max}^t$ \textbf{then break}.\; 
				Set the stepsize $ \eta_{t,i}\gets\min\!\left\{\frac{a_{t,i}}{2L},
				\frac{E_{\max}^t-E_{t,i}}{2L\norm{\nabla f(x_{t,i})}^2}\right\}$ and update $x_{t,i+1}\gets x_{t,i}-\eta_{t,i} \nabla f(x_{t,i})$.\; 
				Set $E_{t,i+1}\gets E_{t,i}+2L\eta_{t,i}\norm{\nabla f(x_{t,i})}^2$ and  $i\gets i+1$.\;
			}
			Set $\widehat{x}_{t+1} = x_{t,i}$ and $t\gets t+1$. 
		}
		\Return $\widehat{x}_{t}$
	\end{algorithm}

	Though the algorithm appears a little bit complicated, its update essentially takes the simple gradient descent form $x\leftarrow x-\eta\nabla f(x)$. The algorithm works as an adaptive mechanism to select the stepsize $\eta$. Its presentation keeps multiple evaluations of the gradient and its norm in one step for clarity, yet with appropriate implementation, the method only needs to compute and save one iterate $x_{t,i}$ and gradient $\nabla f(x_{t,i})$ in memory, all other quantities are scalar values that can be cheaply computed, reused, and updated. 
	
	Compared to existing methods that directly use the Chebyshev schedule in convex quadratic optimization, a major novelty of our method mainly lies in the two levels of adaptivity in the stepsize selection. At the beginning of the $t$-th epoch, we select the largest $m\in\{4,8,16,\ldots\}$ such that 
	\begin{equation}
		\label{eq:nonconvex-cert}
		{MG_th_m^2p^2_{m,\pi/2}} \leq {L^2}/{3},
	\end{equation}
	where $G_t = \norm{\nabla f(\widehat{x}_t)}$ is defined in the algorithm,  $h_m$ and $p_{m,\pi/2}$ are given in Section \ref{sec:prelim}.  This step determines the potential aggressiveness of the Chebyshev schedule in this epoch, because potentially the largest stepsize in this epoch is $\frac{p_{m,\pi/2}}{2L} =\Theta(m^2/2L)$. Then a smaller starting gradient $G_t$ allows the algorithm to use a larger $m$ with more aggressive $\Acal_{m,\pi/2}$.  If no such $m$ exists, we will take the standard GD step with $\eta = 3L/2$, where the factor $3/2$ can be replaced by any number below 2.

	But if we do get some valid $m\geq 4$, then our second level of adaptivity starts to adjust the stepsize. Denote $\Acal_{m,\pi/2} = [a_{t,0},\ldots,a_{t,m-1}]$, and initialize the update budget by
	\begin{equation}
		\label{eq:Ecap}
		E_{t,0}:=0,\qquad
		E_{t,i}:=2L\sum_{j=0}^{i-1}\eta_{t,j}\norm{\nabla f(x_{t,j})}^2,
		\qquad
		E_{\max}^t:=\frac{\theta^2G_t^2p_{m,\pi/2}}{2}.
	\end{equation}
	Then if the previous steps have already made sufficient progress so that $\norm{\nabla f(x_{t,i})}<\theta G_t$, then the current epoch is terminated and we will reinitialize the next epoch. Otherwise, we will select the actual stepsize at the $i$-th update in the following clipped form 
	\begin{equation}
		\label{eq:executed-eta}
		\eta_{t,i} = \min\!\left\{\frac{a_{t,i}}{2L}, \frac{E_{\max}^t-E_{t,i}}{2L\norm{\nabla f(x_{t,i})}^2}\right\}
	\end{equation}
	and then increase $E_{t,i}$ after the gradient descent update of the iterate. That is, even after a Chebyshev schedule $\Acal_{m,\pi/2}$ has been adaptively determined, the actual stepsizes still need to be adaptively controlled. If the previous steps are already aggressive enough so that the budget $E_{\max}^t$ is insufficient, then a planned Chebyshev step is clipped to the remaining budget. In this case the clipped update satisfies $E_{t,i+1}=E_{\max}^t$, and the epoch terminates at its output. This discussion also reveals that the clipping of \eqref{eq:executed-eta} may possibly happen only at the last update of this epoch, while the previous steps take the full Chebyshev step with $\eta_{t,i} = a_{t,i}/2L$.

	These two levels of adaptivity determine how aggressive an epoch of the Chebyshev schedule can potentially be and how much of this schedule can be safely executed, which forms the main algorithmic difference from existing methods that use predetermined Chebyshev schedules.

	\subsection{The Complexity of NACC}
	In this subsection, we establish the complexity of NACC for finding an $\epsilon$-stationary point. We start by establishing the sufficient decrease for \emph{every iteration}  of our method, which is a significant different feature compared to the existing non-monotonically descent accelerations, such as the CUPG method in \citet{carmon2017convex}, as well as the restarted AGD and restarted heavy ball in \citet{li2023restarted}. In particular, we present the main optimization results in the main paper, while relegating the more involved technical analysis of the Chebyshev schedule and Hessian variation to the appendix.
	
	To simplify the notation in this part, we apply the notation of $p_m = p_{m,\pi/2}$ whenever $\phi=\pi/2$ is adopted in the discussion. And we will use $g_{t,i} = \nabla f(x_{t,i})$ as the shorthand of the gradient.

	\begin{lemma}
		\label{lem:nonconvex-path}
		Fix an epoch $t$ with a valid $m\geq 4$, and suppose this epoch terminates within $1\leq k\leq m$ steps. Given $G_t=\norm{\nabla f(x_{t,0})}$, then   
		\begin{equation}
			\sum_{i=0}^{k-1}\norm{x_{t,i+1}-x_{t,i}}
			=\sum_{i=0}^{k-1}\eta_{t,i}\norm{\nabla f(x_{t,i})}
			\leq\frac{\theta p_mG_t}{4L},\nonumber
		\end{equation}
		and we have $\sum_{i=0}^{k-1}\eta_{t,i}\leq {p_m}/{4L}$.
	\end{lemma} 
	\begin{proof} 
		By the definitions in \eqref{eq:Ecap}, the clipping rule \eqref{eq:executed-eta}, and the Line 10-12 of Algorithm \ref{alg:nonconvex}, we know 
		the all $k$ iterations of the epoch satisfy $$E_{t,k}=2L\sum_{i=0}^{k-1}\eta_{t,i}\norm{g_{t,i}}^2\leq E_{\max}^t=\frac{\theta^2p_mG_t^2}{2}.$$
		Also note that for every non-terminating iterate of the epoch $x_{t,i}$ with $0\leq i\leq k-1$, Line 10 of Algorithm \ref{alg:nonconvex} implies $\norm{\nabla f(x_{t,i})}>\theta G_t$.
		Combining this with the above inequality gives 
		$$\sum_{i=0}^{k-1}\eta_{t,i}
		\leq\frac{\sum_{i=0}^{k-1}\eta_{t,i}\norm{g_{t,i}}^2}{\theta^2G_t^2}=\frac{E_{t,k}}{2L\theta^2G_t^2}\leq\frac{E_{\max}^t}{2L\theta^2G_t^2}=\frac{p_m}{4L},$$
		$$\sum_{i=0}^{k-1}\eta_{t,i} \norm{g_{t,i}} \leq \frac{\sum_{i=0}^{k-1}\eta_{t,i}\norm{g_{t,i}}^2}{\theta G_t} = \frac{E_{t,k}}{2L\theta G_t} \leq\frac{E_{\max}^t}{2L\theta G_t} =\frac{\theta p_mG_t}{4L},$$ 
		which completes the proof of this lemma. 
	\end{proof} 
	
	Ignoring the fact that Lemma \ref{lem:nonconvex-path} only provides upper bounds, informally speaking, it characterizes the total stepsize and moving distance within one epoch, which both scales proportional to $p_m = \Theta(m^2)$. For an epoch with at most $m$ steps, this is already $m$-times more aggressive than standard constant GD stepsize. As $m$ is determined by \eqref{eq:nonconvex-cert}, we know $m=\widetilde{\Omega}(L^{1/2}M^{-1/4}G_t^{-1/4})$, which informally explains the source of the $\widetilde{O}(L^{1/2}M^{-1/4}\epsilon^{-1/4})$ times acceleration of our method compared to  GD with $O(1/L)$ stepsizes.

	The next lemma is the key step of the analysis, which guarantees the stability of our methods. Despite the simple form, its proof is highly nontrivial and is hence relegated to the following Subsection \ref{sec:positive-curvature-control}.

	\begin{lemma}
		\label{lem:positive-curvature}
		Suppose Assumption \ref{assumption:f} holds. Fix an epoch $t$ with a valid $m\geq 4$, and suppose this epoch terminates in $1\leq k\leq m$ steps. Then for any $0\leq i\leq k-1$, we have   
		$$\bigg\|\left[\frac{a_{t,i}}{2L}\nabla^2 f(\widehat{x}_t)\right]_+^{1/2}\nabla f(x_{t,i})\bigg\| \leq \frac{33}{25}\norm{\nabla f(x_{t,i})},$$
		where $a_{t,i}$ is $(i+1)$-th entry of the Chebyshev schedule $\Acal_{m,\pi/2}$, see Line 7 of Algorithm \ref{alg:nonconvex}. 
	\end{lemma}
	Together with Lemma \ref{lem:nonconvex-path}, we obtain the following per-iterate sufficient descent result.

	\begin{lemma}
		\label{lem:nonconvex-descent}
		Given Assumption \ref{assumption:f}, for any epoch $t$ with $m\geq 4$, all of its iterates satisfy
		\begin{equation}
			f(x_{t,i+1})\leq f(x_{t,i})-\frac{\eta_{t,i}}{8}\norm{\nabla f(x_{t,i})}^2, \nonumber
		\end{equation}
		where $\eta_{t,i} = \min\!\left\{\frac{a_{t,i}}{2L}, \frac{E_{\max}^t-E_{t,i}}{2L\norm{\nabla f(x_{t,i})}^2} \right\}$ is defined by \eqref{eq:executed-eta}. 
	\end{lemma}
	\begin{proof}
		For an update $x_{t,i+1}=x_{t,i}-\eta_{t,i}\nabla f(x_{t,i})$, let us adopt the shorthands $g_{t,i}=\nabla f(x_{t,i})$, $\widehat{H}_t = \nabla^2 f(\widehat{x}_t)$, and $H_{t,i} = \nabla^2f(x_{t,i})$ in this proof for simplicity. Then we have 
		\begin{equation}
			\label{eq:anchored-third-order}
			\begin{aligned}
				f(x_{t,i+1}) & \overset{\mathrm{(i)}}{\leq} f(x_{t,i}) - \eta_{t,i} \norm{g_{t,i}}^2 + \frac{\eta_{t,i}^2}{2}\big(g_{t,i}\big)^{ \top} \Big(\widehat{H}_t -\widehat{H}_t +H_{t,i}\Big) g_{t,i}  +\frac{M\eta_{t,i}^3}{6} \norm{g_{t,i}}^3\\
				&\overset{\mathrm{(ii)}}{\leq} f(x_{t,i})-\eta_{t,i}\norm{g_{t,i}}^2 +\frac{\eta_{t,i}^2}{2}(g_{t,i})^\top\big[\widehat{H}_t\big]_+g_{t,i} + \frac{M\eta_{t,i}^2}{2}\bigl(\norm{x_{t,i} - \widehat{x}_t} + \eta_{t,i}\norm{g_{t,i}}\bigr)\norm{g_{t,i}}^2\\
				&\overset{\mathrm{(iii)}}{\leq}
				f(x_{t,i})-\eta_{t,i}\bigg(
				1-\frac12\left(\frac{33}{25}\right)^2
				-\frac{\theta MG_tp_m^2}{16L^2}\bigg)\norm{g_{t,i}}^2\\
				&\overset{\mathrm{(iv)}}{\leq}
				f(x_{t,i})-\frac{\eta_{t,i}}{8}\norm{g_{t,i}}^2.
			\end{aligned}
		\end{equation}
		Here both (i) and (ii) use the assumption that $\nabla^2 f$ is $M$-Lipschitz continuous. In particular, (i) is due to the cubic descent lemma  \eqref{eq:cubic} and (ii) directly uses the Lipschitz continuity and the fact that $\widehat{H}_t\preceq\big[\widehat{H}_t\big]_+$.  In (iii), we bound the quadratic term by Lemma~\ref{lem:positive-curvature} and $\eta_{t,i}\leq \frac{a_{t,i}}{2L}$, and we use the triangle inequality,  Lemma~\ref{lem:nonconvex-path}, and $\eta_{t,i}\leq \frac{p_m}{2L}$ to bound 
		$$\frac{M\eta_{t,i}}{2}\Big(\norm{x_{t,i}-\widehat{x}_t}+\eta_{t,i}\norm{g_{t,i}}\Big) \leq \frac{Mp_m}{4L} \sum_{j=0}^i\eta_{t,j}\norm{g_{t,j}} \leq  \frac{\theta MG_tp_m^2}{16L^2}.$$
		Finally, (iv) is due to $h_m\geq3$ for all $m\geq 4$ and the inequality \eqref{eq:nonconvex-cert}.   
	\end{proof}
	
	This lemma shows that for all the epochs that execute the clipped Chebyshev schedule with some $m\geq4$, its iterates monotonically decrease the objective function, and satisfy the standard sufficient descent property even when the stepsize is much larger than $1/L$.  On the other hand, if an epoch does not have a valid $m\geq4$ and executes the constant stepsize GD in Line 5 of Algorithm \ref{alg:nonconvex}, the descent lemma of $L$-smoothness gives 
	\begin{equation}
		\label{eq:fallback-descent}
		f(\widehat{x}_{t+1}) = f\left(\widehat{x}_{t} - \frac{3}{2L}\nabla f(\widehat{x}_{t})\right) \leq f(\widehat{x}_{t}) - \frac{3}{8L}\norm{\nabla f(\widehat{x}_{t})}^2.
	\end{equation}
	Hence, the full algorithm monotonically decreases the objective value. Next, for the epochs that execute the Chebyshev schedule, combining Lemma \ref{lem:nonconvex-descent} and the design of Line 7-12 of Algorithm \ref{alg:nonconvex}, we obtain the following lemma, which shows that these epochs either achieve a $\theta$ times contraction in gradient size or a sufficiently large descent of order $p_mG_t^2/L$.  
	
	\begin{lemma}\label{lem:nonconvex-block}
		Given Assumption~\ref{assumption:f}, for any epoch $t$ with a valid $m\geq4$, unless Algorithm \ref{alg:nonconvex} terminates with $\norm{\nabla f(x)}\leq\epsilon$ within an epoch, its output $\widehat{x}_{t+1}$ satisfies at least one of
		\begin{equation}
			\norm{\nabla f(\widehat{x}_{t+1})} \leq\theta\norm{\nabla f(\widehat{x}_t)} 
			\qquad\text{or}\qquad f(\widehat{x}_t)-f(\widehat{x}_{t+1})\geq\frac{\theta^2p_m}{32L}\norm{\nabla f(\widehat{x}_t)}^2.\nonumber
		\end{equation} 
	\end{lemma}
	\begin{proof}
		Suppose that the epoch runs for $1\leq k\leq m$ steps with $x_{t,0}=\widehat{x}_t$ and $x_{t,k}=\widehat{x}_{t+1}$. If the output $x_{t,k}$ is generated by activating the epoch's stopping criteria $\norm{\nabla f(x_{t,k})} \leq \theta G_t$ in Line 10 of Algorithm \ref{alg:nonconvex}, then the lemma holds with the first inequality being valid. Otherwise, summing up Lemma~\ref{lem:nonconvex-descent} and using the definition of $E_{t,k}$ in \eqref{eq:Ecap}, we obtain
		$$
		f(\widehat{x}_t)-f(\widehat{x}_{t+1})
		\geq\frac18\sum_{i=0}^{k-1}\eta_{t,i}\norm{g_{t,i}}^2
		=\frac{E_{t,k}}{16L}.
		$$
		If $E_{t,k}=E_{\max}^t$, the second inequality holds by  $E_{\max}^t=\theta^2p_mG_t^2/2$, see \eqref{eq:Ecap}. If $E_{t,k}<E_{\max}^t$, then we must have $k=m$. No update has been clipped $\eta_{t,i}=\frac{a_{t,i}}{2L}$ and all gradients are sufficiently large $\norm{g_{t,i}}>\theta G_t$, for all $0\leq i\leq m-1$. As a result, we can obtain  
		$$E_{t,m} = \sum_{i=0}^{m- 1} a_{t,i} \norm{g_{t,i}}^2 > \theta^2G_t^2 \sum_{i=0}^{m-1} a_{t,i} = 2\theta^2m^2G_t^2
		\geq\theta^2p_mG_t^2 = 2E_{\max}^t,$$
		where we used \eqref{eq:general-mass} and $p_m\leq2m^2$. This result contradicts the fact that $E_{t,k}\leq E_{\max}^t$ and the scenario $E_{t,k}<E_{\max}^t$ is impossible. Hence we complete the proof.  
	\end{proof}
	
	Combining the above lemmas, we observe a very clear descent or contract structure in the algorithm's progress. Denote $\Delta_f:=f(\widehat{x}_0)-f^*$, the following Theorem \ref{thm:nonconvex-main} states the final complexity result.
	
	\begin{theorem}\label{thm:nonconvex-main}
		Suppose that Assumption \ref{assumption:f} hold.  Algorithm~\ref{alg:nonconvex} can output an $\epsilon$-stationary point of the problem with at most $$O\left(\frac{\Delta_fM^2}{L^3}
		+ \frac{L^{1/2}}{M^{1/4} \epsilon^{1/4}}\ln\left(\frac{L^2}{M\epsilon}\right) + \frac{\Delta_f L^{1/2} M^{1/4}} {\epsilon^{7/4}}
		\sqrt{\log_2\left(1+\frac{4L^2}{M\epsilon}\right)} \right)$$
		gradient descent steps, where $\widetilde{O}(\epsilon^{-7/4})$ is the leading term.  
	\end{theorem}
	
	\begin{proof}
		Suppose Algorithm~\ref{alg:nonconvex} terminates in epoch $K$. We first deal with the completed epochs $t=0,\ldots,K-1$ before termination, which are partitioned into three sets: 
		\begin{align*}
			\mathcal K_1 &:= \left\{t:\text{ epoch $t$ that executes Line 5}
			\right\},\\
			\mathcal K_2 &:= \left\{t:\text{ epoch $t$ that executes Line 10 by triggering } E_{t,i}=E_{\max}^t\right\},\\
			\mathcal K_3 &:= \left\{t\notin\mathcal{K}_2 : \text{epoch $t$ that executes Line 10 by triggering } \norm{\nabla f(x_{t,i})} \leq \theta G_t\right\}.
		\end{align*}
		For epoch $t\in\mathcal{K}_2\cup\mathcal{K}_3$, we denote their $m$ by $m_t$ for clarity. Next, we bound the number of iterations in each of the three classes. After that, we deal with the final epoch $K$ separately. \vspace{0.2cm}
		
		\noindent\textbf{Epochs in $\mathcal K_1$.}\,\, For $t\in\mathcal K_1$,  no $m\geq 4$ satisfies
		\eqref{eq:nonconvex-cert}. Therefore,  $\norm{\nabla f(\widehat{x}_t)}\geq\frac{L^2}{3Mp_4^2h_4^2}$. By \eqref{eq:fallback-descent}, $f(\widehat{x}_t) - f(\widehat{x}_{t+1}) \geq \frac{L^3}{24M^2p_4^4h_4^4} = \Omega(1)$, which further implies    
		\begin{equation}
			\label{eq:K1-bound}
			|\mathcal K_1| \leq \frac{24\Delta_fM^2p_4^4h_4^4}{L^3} = O\left(\frac{M^2\Delta_f}{L^3}\right).
		\end{equation} 
		
		\noindent\textbf{Epochs in $\mathcal K_2$.}\,\,
		By the way $m_t$ is computed, we know $m_t$ satisfies \eqref{eq:nonconvex-cert} while $2m_t$ not: 
		$$p_{m_t}^2h_{m_t}^2 \leq \frac{L^2}{3MG_t} < p_{2m_t}^2h_{2m_t}^2 \leq 28p_{m_t}^2h_{m_t}^2,$$
		where we used the fact that $\frac{p_{2m_t}^2h_{2m_t}^2}{p_{m_t}^2h_{m_t}^2}\leq 28$ for $m\geq4$. As $\frac{p_m}{m^2}\in[1.6,1.65]$ and $h_m\geq 3$ for $m\geq4$, the left hand side of the above inequality gives $\frac{L^2}{3MG_t}\geq3^2\times1.6^2 m_t^4$ and consequently
		$$\log_2\left(\frac{L^2}{64MG_t}\right)\geq 4\log_2(m_t) = 4 (h_{m_t}-1).$$
		Rearranging the terms gives $4h_{m_t}\leq \log_2\big(\frac{4L^2}{MG_t}\big)$. Then substituting this bound into the right-hand side of the previous inequality gives 
		\begin{equation}\label{eq:nonconvex-m-lower} \frac{\bigl(L^2/MG_t\bigr)^{1/4}} {\sqrt{\log_2\!\left(4L^2/MG_t\right)}}
			\leq 2m_t.
		\end{equation}  
		Also note that for $t\in\mathcal{K}_2$, denote $\delta_t:=f(\widehat{x}_t)-f(\widehat{x}_{t+1})$ for simplicity, the analysis of Lemma \ref{lem:nonconvex-block} shows that epoch terminates with $E_{t,i} = E_{\max}^t$ satisfies $\delta_t \geq \frac{\theta^2p_{m_t}G_t^2}{32L}$. Combined with the above lower bound of $m_t$, we obtain the following finer characterization of it:
		$$m_t\leq \frac{p_{m_t}}{1.6m_t} \leq \frac{20L\delta_t}{\theta^2G_t^2m_t}\leq \frac{40}{\theta^2}G_t^{-7/4}L^{1/2}M^{1/4}\sqrt{\log_2\!\left(4L^2/MG_t\right)}\cdot\delta_t.$$
		As $G_t>\epsilon$ for all $t\in\mathcal{K}_2$ and $\sum_{t\in\mathcal{K}_2}\delta_t\leq \Delta_f$, summing up this inequality over $t\in\mathcal K_2$ gives 
		\begin{equation}
			\label{eq:K2-total}
			\sum_{t\in\mathcal K_2}m_t \leq O\left(\Delta_fL^{1/2}M^{1/4}\epsilon^{-7/4}
			\sqrt{\log_2\left(1+\frac{4L^2}{M\epsilon}\right)}\right).
		\end{equation}
		
		\noindent\textbf{Epochs in $\mathcal K_3$.}
		For the epoch set $\mathcal K_3$, we can partition it into a number of consecutive blocks. Now let us focus on an arbitrary block, and let $t$ be its last epoch, then its $m_t$ satisfies \eqref{eq:nonconvex-cert}, that is $MG_tp_{m_t}^2h_{m_t}^2\leq L^2/3$. Because $G_{t+1}\leq\theta G_t \leq G_t$, we know the next epoch will satisfy $m_{t+1}\geq m_t$ as the current $m_t$ is already feasible to \eqref{eq:nonconvex-cert} with gradient size $G_{t+1}$. Therefore, $t+1\notin\mathcal{K}_1$ must hold. That is, these blocks either ends with the next epoch being the terminating epoch $K$, or they end with the next epoch in $\mathcal{K}_2$. 
		
		For the first case where the block ends at the terminating epoch $K$. The number of this type of block is at most 1. For any epoch $t$ in this block, we use \eqref{eq:nonconvex-cert} to bound $m_t$ by
		\begin{equation} 
			\label{eq:final-block}
			m_t \leq \left(p_{m_t}^2\right)^{1/4} \leq \left(\frac{L^2}{3Mh_{m_t}^2G_t}\right)^{1/4} \leq O\left(\frac{L^{1/2}}{(M\epsilon)^{1/4}}\right).  
		\end{equation}
		Together with the fact that starting epoch $t_0$ of this block should satisfy $G_{t_0} \leq \frac{L^2}{3Mp_4^2h_4^2}$ as it has a valid $m_{t_0}\geq4$. Therefore, this block can contain at most $\mathcal{O}\big(\ln\frac{L^2}{M\epsilon}\big)$ epochs, and the total number of gradient steps contained in this block is at most $\widetilde{O}(\epsilon^{-1/4})$.
		
		Then, for the second case where the block ends at an epoch that belongs to $\mathcal{K}_2$, let us focus on an arbitrary block of this type. For any $t$ in this block, as $\frac{p_{2m}^2h_{2m}^2}{p_{m}^2h_{m}^2}\leq 28$ for $m\geq4$, there is an absolution constant integer $n_0>0$ such that $28\theta^{n_0}\leq 1$. If this block does not end in $n_0$, then the $m$ parameter will double after at most $n_0$ steps because  
		$$MG_{t+n_0}p_{2m_t}^2h_{2m_t}^2\leq 28\theta^{n_0}\cdot MG_{t}p_{m_t}^2h_{m_t}^2\leq L^2/3$$
		will satisfy \eqref{eq:nonconvex-cert}.
		That is, $m_{t+n_0}\geq2m_t$. Now suppose this block ends with some $t'\in\mathcal{K}_2$, and we group the epochs in this block backwards into sub-blocks of at most $n_0$ epochs. Then their $m$'s can be bounded by
		$m_{t'},\frac{m_{t'}}{2}, \frac{m_{t'}}{2^2},\ldots.$ Therefore, the total number of steps in this block is at most
		\begin{equation} 
			n_0 \cdot \sum_{j=0}^{\infty}2^{-j}m_{t'}
			= 2n_0m_{t'}. \nonumber
		\end{equation} 
		Following the above discussion, the total number of gradient steps in this type of blocks can be upper bounded by  
		\begin{equation} 
			2n_0\sum_{t'\in\mathcal K_2}m_{t'} \leq O\left( \Delta_f L^{1/2} M^{1/4} \epsilon^{-7/4}
			\sqrt{\log_2\left(1+\frac{4L^2}{M\epsilon}\right)}\right), \nonumber
		\end{equation}
		where we use the inequality \eqref{eq:K2-total}. Combining the bounds for the two types of blocks upper bounds the total number of steps for epochs in $\mathcal{K}_3$. \vspace{0.2cm}
		
		\noindent\textbf{Final epoch $K$.}\,\, For $t=K$, the inequality \eqref{eq:final-block} still holds. \vspace{0.2cm}
		
		\noindent Taken together, the upper bounds established in the four cases above prove the theorem.
	\end{proof}
	
	\subsection{The Key Positive Curvature Bound of Lemma \ref{lem:positive-curvature}}
	\label{sec:positive-curvature-control}
	In this subsection, we provide the key analysis of Lemma \ref{lem:positive-curvature}, which was skipped previously due to the heavy technicality. Throughout this section, we will work within one single epoch. Therefore, we can drop the subscript $t$ for a simpler presentation, and denote
	$$x_j=x_{t,j},\quad x_0=\widehat{x}_t,\quad a_j=a_{t,j},\quad \eta_j=\eta_{t,j}, \quad g_j=\nabla f(x_{t,j}),\quad G=\norm{\nabla f(\widehat{x}_t)}.$$
	Let us inherit the setting of Lemma \ref{lem:positive-curvature} so that this epoch contains $0\leq k\leq m$ steps. Consider a specific iterate $x_i$ for some $0\leq i \leq k-1$ in this epoch. By the analysis of Lemma \ref{lem:nonconvex-block}, all previous iterates are unclipped as clipping may only (possibly) happen at the end of each epoch. That is,   $\eta_j=\frac{a_j}{2L}$ for all $0\leq j<i$. 
	
	Because $x_{j+1} = x_j - \eta_j \nabla f(x_{j})$ is generated by the gradient descent step, we have 
	$$\nabla f(x_{j+1})-\nabla f(x_{j}) = -\int_0^1\eta_j\nabla^2f\bigl(x_j-s\eta_j\nabla f(x_{j})\bigr)\nabla f(x_{j})\dd s.$$ 
	Denote the Hessian and its averaged variants as
	\begin{equation}
		\widehat H:=\nabla^2f(x_0)  \qquad
		\overline   H_j := \int_0^1 \nabla^2f\bigl(x_j+s(x_{j+1}-x_j)\bigr)\,\dd s, \quad 0\leq j<i. \nonumber
	\end{equation}
	Then the above gradient difference equation can be written as an exact nonautonomous linear dynamical system with time varying matrix  $\overline{H}_j$:
	\begin{equation}\label{eq:epoch-gradient-dynamics}
		g_{j+1}=(I-\eta_j\overline H_j)g_j=(I-\eta_j\widehat H)g_j+\eta_j(\widehat H-\overline H_j)g_j,
		\quad 0\leq j<i.
	\end{equation} 
	For $j\geq1$, define the accumulated Hessian variation by
	\begin{equation}\label{eq:epoch-hessian-variation}
		V_j:=\norm{\overline H_0-\widehat H}+\sum_{j'=1}^{j-1}\norm{\overline H_{j'}-\overline H_{j'-1}},
	\end{equation}
	and default $V_0:=0$ when $j=0$. Define $\widehat{H}_+ = \big[\widehat{H}\big]_+\succeq0$ and $\widehat{H}_- = \widehat{H}_+ -\widehat{H}\succeq0$, then we have $\widehat{H} = \widehat{H}_+-\widehat{H}_-$. Define $\Pi_+$ and $\Pi_-$ the orthogonal projections to the eigenspaces of $\widehat H_+$ and $\widehat H_-$, respectively, and define $g_j^+:=\Pi_+ g_j$ and $g_j^-:=\Pi_- g_j$,
	then the following lemma holds. 
	
	\begin{lemma}
		\label{lem:curvature-weighted-gradient}
		Under the setting of Lemma \ref{lem:positive-curvature}, it holds that 
		\begin{equation} 
			\left\| \sqrt{\frac{a_i}{2L}}\cdot \widehat H_+^{1/2}g_i \right\| \leq\frac76 \cdot \max_{0\leq j\leq i} \norm{g_j^+} + \omega_i  V_i \cdot \max_{0\leq j\leq i}\norm{g_j}, \nonumber
		\end{equation}
		where $\omega_i:=\frac{\pi m^2}{2\sqrt2L}(h_m-1).$
	\end{lemma}
	
	\begin{proof}
		By the recursive binary tree structure of $\Acal_{m,\pi/2}$, see Figure \ref{fig:tree-m4} for example, there exists a maximal sub-tree with $a_i$ being the right-most leaf,
		which we write as
		$$Q=[a_s,\ldots,a_i]=\Acal_{n,\phi}\quad\mbox{with}\quad n = i-s+1\leq m.$$
		Such a maximal sub-tree can by determined by tracing from the leaf $a_i$ upward until it reaches a left-child node. If the leaf $a_i$ itself is already a left-child node, then $s=i$ and this block contains a single element $a_i$. Then that node will be the $\Acal_{n,\phi}$ that we are looking for. According to \citet{lebedev1971ordering,agarwal2021acceleration}, for this block, we must have $\phi\in[\pi/2,\pi)$ and $a_i=p_{n,\phi}$. 
		
		Projecting  
		\eqref{eq:epoch-gradient-dynamics} to the eigenspace of
		$\widehat H_+$ and unrolling from $s$ to $i$ gives
		$$g_i^+ = \Phi_{s,i}g_s^+ + \sum_{j=s}^{i-1} \eta_j \Phi_{j+1,i}\Pi_+e_j\quad\mbox{with}\quad e_j:=(\widehat H-\overline H_j)g_j.$$
		Here, the matrix $\Phi_{j,i} := \prod_{\ell=j}^{i-1} (I-\eta_\ell \widehat H_+)$ for $\forall j<i$, and we default $\Phi_{i,i}:=I$ for $j=i$. As  $\eta_j={a_j}/{2L}$ for every $s\leq j<i$, we have  
		\begin{equation}
			\label{eq:c-w-g-1}
			\left\| \sqrt{\frac{a_i}{2L}} 
			\widehat H_+^{1/2}g_i \right\|
			= \left\| \sqrt{\frac{a_i}{2L}} 
			\widehat H_+^{1/2}g_i^+ \right\| \leq \left\| \sqrt{\frac{a_i}{2L}}  \widehat H_+^{1/2} \Phi_{s,i}g_s^+ \right\| \!+\!
			\sum_{j=s}^{i-1} \frac{a_j}{2L} \left\|
			\sqrt{\frac{a_i}{2L}}  \widehat H_+^{1/2}\Phi_{j+1,i}\Pi_+e_j \right\|\!. 
		\end{equation} 
		By spectral calculus, Lemmas~\ref{lem:preparation-moment} and~\ref{lem:weighted-response} in Appendix~\ref{sec:cheb-toolkit}, applied to $\widehat H_+/(2L)$, give the following two inequalities, respectively.
		\begin{equation}
			\label{eq:c-w-g-2}
			\left\| \sqrt{\frac{a_i}{2L}}\,
			\widehat H_+^{1/2}\Phi_{s,i} \right\| \leq \frac{7}{6}.
		\end{equation}
		\begin{equation}
			\label{eq:c-w-g-3}
			\sum_{j=s}^{i-1} \frac{a_j}{2L}
			\left\| \sqrt{\frac{a_i}{2L}}\,
			\widehat H_+^{1/2}\Phi_{j+1,i} \right\| \leq \frac{\pi n^2}{2\sqrt2L}(h_n-1) \leq \frac{\pi m^2}{2\sqrt2L}(h_m-1).
		\end{equation}
		Substituting \eqref{eq:c-w-g-2}, \eqref{eq:c-w-g-3}, and the bound $\|\Pi_+e_j\|\leq \norm{\overline H_j-\widehat H}\cdot\|g_j\|\leq V_i\cdot \max_{0\leq j\leq i}\norm{g_j}$ for all $j<i$, to the inequality \eqref{eq:c-w-g-1} finishes the proof of this lemma. 
	\end{proof}
	Therefore, it remains to upper bound the Hessian variation term $V_i$, and the maximal $\|g_j^+\|$ and $\|g_j^-\|$ along the iterations of an epoch. 
	
	\begin{lemma} 
		\label{lem:hessian-variation}
		Under the setting of Lemma \ref{lem:positive-curvature}, we have $V_i\leq \frac{\theta Mp_mG}{4L}$ for all $i\leq k$.
	\end{lemma}
	\begin{proof}
		The case $i=0$ is immediate. For $1\leq j < i$, the  definition of $\overline{H}_j$ and the  Hessian $M$-Lipschitz continuity yields 
		\begin{equation*}
			\begin{aligned}
				\norm{\overline H_j-\overline H_{j-1}} & =  \left\| \int_0^1\! \left[ \nabla^2 f\bigl(x_j + s(x_{j+1}-x_j)\bigr) -\nabla^2 f \bigl(x_{j-1}+s(x_j-x_{j-1})\bigr)\right]\dd s \right\| \\
				& \leq  M\int_0^1\big\|s(x_{j+1}-x_j)+(1-s)(x_j-x_{j-1})\big\|\,\dd s \\
				& \leq M\int_0^1s\big\|x_{j+1}-x_j\big\|\,\dd s +    M\int_0^1(1-s)\big\|x_j-x_{j-1}\big\|\,\dd s \\
				& = \frac M2\Big(\norm{x_{j+1}-x_j}+\norm{x_j-x_{j-1}}\Big)
			\end{aligned}
		\end{equation*}
		Through similar analysis, we also have $\norm{\overline H_0-\widehat H} \leq \frac M2\norm{x_1-x_0}.$ Then summing up the above terms of $V_i$ given in \eqref{eq:epoch-hessian-variation} derives 
		\begin{equation}
			\label{eq:hessian-prefix-variation}
			V_i \leq M\sum_{j=0}^{i-2}\eta_j\norm{g_j}
			+\frac M2\eta_{i-1}\norm{g_{i-1}}.
		\end{equation}
		Finally, applying Lemma \ref{lem:nonconvex-path} to the above bound proves this lemma. 
	\end{proof}
	
	\begin{lemma}
		\label{lem:negative-gradient-bound}
		Under the setting of Lemma \ref{lem:positive-curvature},  $\norm{g_j^-} \leq \norm{g_i^-} + \frac{\theta^2Mp_m^2G^2}{32L^2}$ holds for $\forall 0\leq j\leq i$.
	\end{lemma}
	\begin{proof}
		For $i\geq1$, projecting the gradient dynamics \eqref{eq:epoch-gradient-dynamics} to the eigenspace of $\widehat H_-$ gives
		\begin{equation}
			\label{eq:negative-dynamics}
			g_{j+1}^-=(I+\eta_j\widehat H_-)g_j^-+\eta_j\Pi_-e_j
		\end{equation} 
		where $e_j = (\widehat H-\overline H_j)g_j$ is defined in consistence with the proof of Lemma \ref{lem:curvature-weighted-gradient}. By $\widehat H_-\succeq0$, we have $\norm{(I+\eta_j\widehat H_-)^{-1}}\leq1$, and hence $\norm{g_j^-} \leq \norm{g_{j+1}^-} + \eta_j\norm{e_j}.$ Together with the bound 
		$$\norm{e_j} \leq \norm{\overline H_{j}-\widehat H} \cdot \norm{g_j} \leq V_{j+1}\norm{g_j}
		\leq M\norm{g_j}\bigg(\sum_{\ell=0}^{j-1}\eta_\ell\norm{g_\ell}+\frac12\eta_j\norm{g_j}\bigg),$$
		where the last inequality is by \eqref{eq:hessian-prefix-variation}. Consequently, we have 
		\begin{equation} 
			\label{eq:negative-gradient-bound-1}
			\begin{aligned}
				\max_{0\leq j\leq i} \norm{g_j^-} & \leq \norm{g_i^-} + \sum_{j=0}^{i- 1} \eta_j \norm{e_j} \\
				& \leq \norm{g_i^-} + M\sum_{j=0}^{i-1}\eta_j\norm{g_j}\bigg(\sum_{\ell=0}^{j-1}\eta_\ell\norm{g_\ell}+\frac12\eta_j\norm{g_j}\bigg)\\
				&= \norm{g_i^-} +\frac M2 \bigg(\sum_{j=0}^{i- 1} \eta_j \norm{g_j} \bigg)^2 \leq \norm{g_i^-} + \frac{\theta^2Mp_m^2G^2}{32L^2},
			\end{aligned}
		\end{equation}
		where the last inequality follows from Lemma \ref{lem:nonconvex-path}.  
	\end{proof}
	
	Finally, we bound the positive component $g_j^+$ by the following lemma.  
	
	\begin{lemma}
		\label{lem:positive-gradient-bound}
		Under the setting of Lemma \ref{lem:positive-curvature}, we have $\|g_j^+\|\leq A_1 G + A_2 \|g_i^-\|$, with constants
		$$A_1 = e^{\nu} + \frac{\theta^2 e^{\nu+1}\,\nu(2+\lambda)}{96h_m^2} \quad\mbox{and}\quad
		A_2 = \nu e^{\nu+1} (1+\lambda),$$ 
		where $\nu = \frac{\theta m^2}{24p_m} \big(\frac{1}{2h_m} + \frac{1}{3h_m^2}\big)$ and $\lambda = \ln\big(1+ \frac{p_m}{4}\big)$.
	\end{lemma}
	\begin{proof}
		First, projecting \eqref{eq:epoch-gradient-dynamics} to the eigenspace of $\widehat{H}_+$ gives
		$$g_{j+1}^+=\Pi_+(I-\eta_j\overline H_j)(g_j^++g_j^-) =(I-\eta_j\Pi_+\overline H_j\Pi_+)g_j^+
		-\eta_j\Pi_+\overline H_jg_j^-.$$
		Default $\Psi_{j,i}:=I$ for $j=i$ and define for $\forall j<i$ that
		$$\Psi_{j,i} := (I-\eta_{i-1} \Pi_+\overline H_{i-1}\Pi_+) \cdots (I-\eta_j \Pi_+\overline H_j\Pi_+).$$
		We should note that the matrices in this product are not necessarily commutative and hence the order of the product cannot be exchanged. We also denote $w_j=\Pi_+\overline H_jg_j^-$. Then for any $1\leq j\leq i$, unrolling the above recursive formula gives
		\begin{equation}
			\label{eq:positive-gradient-unrolled} 
			g_j^+ = \Psi_{0,j}g_0^+ - \sum_{\ell=0}^{j-1} \eta_\ell \Psi_{\ell+1,j} w_\ell  = \Psi_{0,j}g_0^+ - F_{0,j} w_0 - \sum_{\ell=1}^{j-1} F_{\ell,j} \left( w_\ell - w_{\ell-1} \right), 
		\end{equation}
		where we use the discrete summation by parts trick with $F_{\ell,j}  = \sum_{q=\ell}^{j-1} \eta_q \Psi_{q+1,j}$ and default $F_{j,j}:=0.$ Utilizing the property of the Chebyshev schedule $\Acal_{m,\pi/2} = [a_0, a_1,\cdots,a_{m-1}]$ and the fact that $\eta_j = a_j/2L$, we have 
		\begin{equation}
			\label{eq:positive-prod-res}
			\norm{\Psi_{0,j}}  \leq \ee^{c_mV_i}, \qquad \norm{F_{\ell,j}} \leq c_m \ee^{1+c_mV_i}, \qquad 0\leq \ell \leq j\leq i,
		\end{equation} 
		where $c_m:=\frac{m^2}{2L} \big( \frac{h_m}{2} + \frac13 \big),$ and the proof of these inequalities is given in Section~\ref{sec:positive-prod-res}. Also note that $\Pi_+\overline{H}_\ell\Pi_- = \Pi_+(\overline{H}_\ell-\widehat{H})\Pi_-$, for all $\ell$, we can obtain
		\begin{equation*} 
			\begin{aligned}
				\norm{w_0}  &+ \sum_{\ell=1}^{i-1} \norm{ w_\ell - w_{\ell-1}}  = \norm{\Pi_+(\overline{H}_0-\widehat{H})\Pi_-g_0^-} + \sum_{\ell=1}^{i-1} \norm{ \Pi_+ \overline{H}_\ell \Pi_- g_\ell^- - \Pi_+\overline{H}_{\ell-1}\Pi_-g_{\ell-1}^-}\\ 
				& \leq \norm{\overline{H}_0-\widehat{H}}\cdot\norm{g_0^-} + \sum_{\ell=1}^{i-1} \Big(\norm{\overline{H}_\ell-\overline{H}_{\ell- 1}} \cdot \norm{g_\ell^-} +  \big\| \overline{H}_{\ell-1} - \widehat{H}\big\|\cdot\|g_\ell^--g_{\ell-1}^-\|\Big)\\
				&\leq V_i\bigg(\max_{0\leq \ell\leq i}\norm{g_\ell^-} + \sum_{\ell=0}^{i-1} \norm{g_{\ell+1}^--g_\ell^-} \bigg).
			\end{aligned}
		\end{equation*}
		Substituting the above bounds to \eqref{eq:positive-gradient-unrolled} and using $\norm{g_0^+}\leq G$, we obtain \begin{equation}
			\label{eq:positive-forced-stability} 
			\norm{g_j^+} \leq \ee^{c_m V_i} G+ \ee^{ 1+ c_m V_i} c_m V_i\bigg(\max_{0\leq \ell\leq i}\norm{g_\ell^-} + \sum_{\ell=0}^{i-1} \norm{g_{\ell+1}^--g_\ell^-} \bigg),\quad \forall j\leq i. 
		\end{equation}
		The $\max_{0\leq \ell\leq i}\norm{g_\ell^-}$ term can be controlled by Lema \ref{lem:negative-gradient-bound}. It remains to bound the remaining summation term. Rearranging the equation \eqref{eq:negative-dynamics} gives $g_{\ell+1}^--g_\ell^-= \eta_\ell\widehat H_- g_\ell^-+\eta_\ell\Pi_-e_\ell$, which implies 
		\begin{equation}
			\label{eq:bound-1}
			\sum_{\ell=0}^{i-1}\|g_{\ell+1}^- - g_\ell^-\| \leq \sum_{\ell=0}^{i-1} \eta_\ell\|\widehat H_- g_\ell^-\| + \sum_{\ell=0}^{i-1}\eta_\ell\|e_\ell\|,
		\end{equation}
		for which the second term has already been bounded in \eqref{eq:negative-gradient-bound-1}. To bound the first term, we again use \eqref{eq:negative-dynamics} to obtain 
		$g_\ell^- = (I+\eta_\ell\widehat H_-)^{-1}\big(g_{\ell+1}^--\eta_\ell\Pi_-e_\ell\big).$ Unrolling it gives  
		$$g_\ell^- = \prod_{q=\ell}^{i-1}(I+\eta_q\widehat H_-)^{-1}g_i^-
		-\sum_{r=\ell}^{i-1}\eta_r \prod_{q=\ell}^{r}(I+\eta_q\widehat H_-)^{-1}\Pi_-e_r, \quad 0\leq\ell<i.$$
		Consequently, the first term of \eqref{eq:bound-1} can be upper bounded by 
		\begin{equation}
			\sum_{\ell=0}^{i-1}\!\eta_\ell\norm{\widehat H_-g_\ell^-}  \leq \norm{g_i^-}\sum_{\ell=0}^{i-1}
			\big\|S_{i-1,\ell}\big\| + \sum_{r=0}^{i-1} \eta_r\norm{e_r} \sum_{\ell=0}^{r} \big\|S_{r,\ell}\big\|. \nonumber
		\end{equation} 
		where we denote $S_{r,\ell} = \eta_\ell\widehat H_- \prod_{q=\ell}^{r}(I+\eta_q\widehat H_-)^{-1}$ for $0\leq\ell\leq r<i$ to simplify the formula. Then the matrix norm factors can be controlled as follows. Fix any $0\leq\ell\leq r<i$. Let $\mu\in[0,L]$ be an arbitrary eigenvalue of $\widehat{H}_-$ and let $\mu_{r,\ell}$ be the corresponding eigenvalue of $S_{r,\ell}$ that shares the same eigenvector with $\mu$. Then  
		$$\mu_{r,\ell} = \frac{\eta_\ell\mu} {\prod_{q=\ell}^{r}(1+\eta_q\mu)} \leq \frac{\eta_\ell\mu} {1+\mu\sum_{q=\ell}^{r}\eta_q} \leq \frac{\eta_\ell} {L^{-1}+\sum_{q=\ell}^{r}\eta_q}.$$
		That is, $\|S_{r,\ell}\|\leq {\eta_\ell}/ \big(L^{-1}+\sum_{q=\ell}^{r}\eta_q\big)$. Consequently, we have   
		$$ 
		\sum_{\ell=0}^{r} \big\|S_{r,\ell}\big\| \leq
		\sum_{\ell=0}^{r} \frac{\eta_\ell}{L^{-1}+\sum_{q=\ell}^{r}\eta_q} \leq \sum_{\ell=0}^{r} \int_{\sum_{q=\ell+1}^{r}\eta_q}^{\sum_{q=\ell}^{r}\eta_q} \frac{\dd v}{L^{-1}+v} = \ln\bigg(1+L\sum_{q=0}^{r}\eta_q\bigg) \leq\lambda, 
		$$
		where the last inequality is by Lemma \ref{lem:nonconvex-path}, which provided $\sum_{q=0}^{r}\leq p_m/4L$. As a result,  
		\begin{equation} 
			\label{eq:bound-2}
			\sum_{\ell=0}^{i-1}\eta_\ell\norm{\widehat H_-g_\ell^-} \leq \lambda\bigg( \norm{g_i^-}+\sum_{\ell=0}^{i- 1}\eta_\ell\norm{e_\ell}\bigg).
		\end{equation} 
		Then combining \eqref{eq:positive-forced-stability}, \eqref{eq:bound-1},   \eqref{eq:bound-2}, and \eqref{eq:negative-gradient-bound-1} from the proof of Lemma \ref{lem:negative-gradient-bound}, we obtain
		$$
		\norm{g_j^+} \leq \ee^{c_mV_i}G + \ee^{1+c_mV_i} c_mV_i \bigg( (1+\lambda)\norm{g_i^-}
		+(2+\lambda)\sum_{\ell=0}^{i-1}\eta_\ell\norm{e_\ell}
		\bigg).
		$$  
		Finally, Lemma \ref{lem:hessian-variation} and \eqref{eq:nonconvex-cert} imply
		\begin{equation}\label{eq:positive-gradient-exponent}
			c_mV_i
			\leq
			\frac{\theta m^2Mp_mG}{8L^2}
			\left(\frac{h_m}{2}+\frac13\right)
			\leq
			\frac{\theta m^2}{24p_m}
			\left(\frac{1}{2h_m}+\frac{1}{3h_m^2}\right)
			=\nu.
		\end{equation}
		The inequality \eqref{eq:negative-gradient-bound-1} in the proof of Lemma \ref{lem:negative-gradient-bound} and  \eqref{eq:nonconvex-cert} gives
		$$
		\sum_{\ell=0}^{i-1}\eta_\ell\norm{e_\ell}
		\leq
		\frac{\theta^2Mp_m^2G^2}{32L^2}
		\leq
		\frac{\theta^2}{96h_m^2}G.
		$$
		Substituting these bounds completes the proof of this lemma. 
	\end{proof}
	
	At this moment, Lemma \ref{lem:hessian-variation}--\ref{lem:positive-gradient-bound} have already provided the bounds needed in Lemma \ref{lem:curvature-weighted-gradient}. Note that the epoch's stopping criteria also guarantees that  $G\leq\norm{g_i}/\theta$. Substituting these bounds into Lemma \ref{lem:curvature-weighted-gradient} and using \eqref{eq:nonconvex-cert} to cancel out $G$ whenever it appears in the computation, the proof of Lemma \ref{lem:positive-curvature} reduces to elementary numerical verification, and is omitted for simplicity. Despite the lengthy analysis, the key mechanism behind it is clear. The negative curvature never hurts the stability of the algorithm, while the positive curvature related terms can be carefully controlled together by the reversed fractal schedule and the adaptive stopping and clipping rules within the epoch. This ensures sufficient descent of the algorithm throughout the epoch even when step sizes are much larger than $1/L$. 
	
	\section{A Unified Acceleration under Additional Convexity Condition}
	\label{sec:convex}
	In this section, we show that the NACC method can automatically accelerate the gradient descent update for convex optimization without any modification. In addition to Assumption~\ref{assumption:f}, formally,  we impose the following convexity assumption.
	\begin{assumption}\label{assumption:convex}
		The function $f$ is convex and attains its minimum at some $x^*\in\mathbb R^d$.
	\end{assumption}
	
	In the follwoing discussion of this section, we denote  $R:=\max\{1,\norm{\widehat{x}_0-x^*}\}$ for simplicity. We keep the notation $p_m=p_{m,\pi/2}$ and $g_{t,i}=\nabla f(x_{t,i})$ in the proof for convenience. In the formal statement of lemmas and theorems, we still stick to $p_{m,\pi/2}$ and $\nabla f(x_{t,i})$.
	
	As the convexity assumption provides us extra structure for the analysis, it facilitates tighter convergence results. Next, we present the key steps for establishing an $\widetilde{O}(\epsilon^{-3/4})$ complexity for finding an approximate solution $x$ such that $f(x)-f(x^*)\leq R\epsilon$, with a few complicated pure technical lemmas relegated to the appendix. In particular, we use the $(R\epsilon)$-optimal solution instead of an $\epsilon$-optimal solution as we want to keep the implementable termination criteria $\norm{\nabla f(\widehat{x}_t)}\leq \epsilon$ of Algorithm \ref{alg:nonconvex}.

	\begin{lemma}\label{lem:convex-local-safety}
		Suppose Assumptions~\ref{assumption:f} and~\ref{assumption:convex} hold. Fix an epoch $t$ of Algorithm~\ref{alg:nonconvex} with a valid $m\geq4$, and suppose this epoch terminates in $1\leq k\leq m$ steps. For every $0\leq i\leq k-1$, if  $\eta_{t,i}>1/L$, then we have
		\begin{equation}
			\eta_{t,i}\nabla f(x_{t,i})^\top\nabla^2f(x_{t,i})\nabla f(x_{t,i})
			+\frac{M\eta_{t,i}^2}{3}\norm{\nabla f(x_{t,i})}^3
			\leq\norm{\nabla f(x_{t,i})}^2,
			\nonumber\end{equation}
		where $\eta_{t,i}$ is the actual stepsize defined by \eqref{eq:executed-eta}.
	\end{lemma}
	
	Together with convexity, this lemma gives the following sufficient descent result, and we can further show that the distance to $x^*$ will not increase after every iterations.
	
	\begin{proposition}\label{prop:convex-descent}
		Given Assumptions~\ref{assumption:f} and~\ref{assumption:convex}, for any epoch $t$ of Algorithm~\ref{alg:nonconvex} that has a valid $m\geq4$, all of its iterates satisfy
		\begin{equation}
			f(x_{t,i+1})\leq f(x_{t,i})-\frac{\eta_{t,i}}2\norm{\nabla f(x_{t,i})}^2.
			\nonumber\end{equation}
		Moreover, every iterate satisfies $\norm{x_{t,i+1}-x^*}\leq\norm{x_{t,i}-x^*}$. 
	\end{proposition}
	\begin{proof}
		For iterates with $\eta_{t,i}\leq1/L$, the standard descent lemma under $L$-smoothness gives $f(x_{t,i+1}) \leq f(x_{t,i})-\frac{\eta_{t,i}}2\norm{g_{t,i}}^2.$ 
		For iterates with $\eta_{t,i}>1/L$, the  inequality \eqref{eq:cubic} gives
		\begin{align*}
			f(x_{t,i+1}) \!\leq\! f(x_{t,i})\!-\!\eta_{t,i}\norm{g_{t,i}}^2
			\!+\!\frac{\eta_{t,i}^2}{2}g_{t,i}^\top\nabla^2\!f(x_{t,i})g_{t,i}
			\!+\!\frac{M\eta_{t,i}^3}{6}\norm{g_{t,i}}^3 \leq f(x_{t,i})\!-\! \frac{\eta_{t,i}}2\norm{g_{t,i}}^2,
		\end{align*}
		where the last inequality is due to Lemma~\ref{lem:convex-local-safety}. This proves the sufficient descent part of the proposition. To show that the distance to $x^*$ does not increase, we expand the square  
		$$
		\norm{x_{t,i+1}-x^*}^2
		=\norm{x_{t,i}-x^*}^2-2\eta_{t,i}\ip{g_{t,i}}{x_{t,i}-x^*}
		+\eta_{t,i}^2\norm{g_{t,i}}^2.
		$$
		By convexity and the first inequality, we have
		$$
		\ip{g_{t,i}}{x_{t,i}-x^*}
		\geq f(x_{t,i})-f^*
		\geq f(x_{t,i})-f(x_{t,i+1})
		\geq\frac{\eta_{t,i}}2\norm{g_{t,i}}^2.
		$$
		Substituting this inequality into the above equation of squared distances proves the result.
	\end{proof}
	
	If an epoch does not have a valid $m\geq4$ and executes the constant stepsize GD in Line 5 of Algorithm \ref{alg:nonconvex}, the sufficient descent in the nonconvex setting \eqref{eq:fallback-descent} still holds. As $3/2L<2/L$, the standard analysis of gradient descent for smooth convex functions \citep{nesterov2013introductory} gives $\norm{\widehat{x}_{t+1}-x^*}\leq\norm{\widehat{x}_t-x^*}$. Together with Proposition~\ref{prop:convex-descent}, this shows that the objective value and the distance to $x^*$ are nonincreasing throughout the algorithm. 
	
	\begin{proposition}\label{prop:convex-block}
		Given Assumptions~\ref{assumption:f} and~\ref{assumption:convex}, for any epoch $t$ of Algorithm~\ref{alg:nonconvex} with a valid $m\geq4$, unless it terminates with $\norm{\nabla f(x)}\leq\epsilon$, its output satisfies at least one of
		\begin{equation}
			\norm{\nabla f(\widehat{x}_{t+1})}\leq\theta G_t,
			\qquad\text{or}\qquad
			f(\widehat{x}_t)-f(\widehat{x}_{t+1})\geq\frac{\theta^2p_mG_t^2}{8L}.
			\nonumber\end{equation}
		The second inequality holds whenever the epoch exhausts its budget, i.e., when $E_{t,k}=E_{\max}^t$, where $k$ is the number of updates executed in epoch $t$.
	\end{proposition}
	As the proof of this proposition follows exactly the same argument of Lemma~\ref{lem:nonconvex-block}, with Lemma~\ref{lem:nonconvex-descent} replaced by Proposition~\ref{prop:convex-descent}. Therefore, we omit its proof.  At this moment, the descent or contract structure of the above result is exactly the same as the nonconvex setting. Therefore, following the same analysis strategy while adapting to the convex prblem, we obtain the following theorem. 
	
	\begin{theorem}\label{thm:convex-main}
		Suppose Assumptions~\ref{assumption:f} and~\ref{assumption:convex} hold. For $\forall \epsilon>0$, Algorithm~\ref{alg:nonconvex} terminates in
		\begin{equation}
			O\!\left( 1 + L^{1/2} M^{1/4} R \epsilon^{-3/4} \sqrt{1+\log_2\!\left(1+\frac{L^2}{M\epsilon}\right)}+ \sqrt{\frac{LR}{\epsilon}} \right) = \widetilde{O} \big(\epsilon^{-3/4}\big)
			\nonumber\end{equation}
		gradient descent steps. Its output point $\bar{x}$ satisfies $\|\nabla f(\bar{x})\|\leq \epsilon$ and $f(\bar{x})-f^*\leq R\epsilon$.
	\end{theorem}
	\begin{proof}
		By proposition~\ref{prop:convex-descent} and the following discussion on constant $3/2L$ stepsize steps,  we know that $\norm{x_{t,i}-x^*}\leq R$ for all the iterations.. Consequently, $\norm{\nabla f(x_{t,i})}\leq LR$ and 
		$$
		f(x)-f^*\leq\ip{\nabla f(x_{t,i})}{x_{t,i}-x^*}\leq R\norm{\nabla f(x_{t,i})}.
		$$
		That is, terminating at some point $\norm{\nabla f(x_{t,i})}\leq\epsilon$ implies $f(x_{t,i})-f^*\leq R\epsilon$. 
		
		Similar to the proof of Theorem \ref{thm:nonconvex-main}, suppose Algorithm 1 terminates in epoch $K$. We first deal with the completed epochs $t=0,\ldots,K-1$ before termination. We partition these epochs into three sets $\mathcal K_1$, $\mathcal K_2$, and $\mathcal K_3$ in the same way as in the proof of Theorem~\ref{thm:nonconvex-main}, and use the same notation $m_t$ for the selected lengths. Write $\Delta_t:=f(\widehat{x}_t)-f^*$ and $\delta_t:=f(\widehat{x}_t)-f(\widehat{x}_{t+1})$. For the ease of notation, in this proof, we will also use $C>0$ denotes an absolute constant that may vary between inequalities. 
		
		\noindent\textbf{Epochs in $\mathcal K_1\cup\mathcal K_2$.}\,\,
		For $t\in\mathcal K_1$, inequality \eqref{eq:fallback-descent} gives $\delta_t\geq\frac{3G_t^2}{8L}$, while the infeasibility of $m=4$ in the criteria \eqref{eq:nonconvex-cert} gives $\frac{MG_t}{L^2}>\frac{1}{3p_4^2h_4^2}$. As a result, we have 
		$$
		G_t^{-7/4}\delta_t \geq G_t^{-7/4}\frac{3G_t^2}{8L} = \frac{3G_t^{1/4}}{8L} > \frac{3}{8L}\left(\frac{L^2}{3Mp_4^2h_4^2}\right)^{1/4} \geq  \frac{C}{L^{1/2}M^{1/4}}.
		$$
		For $t\in\mathcal K_2$, we should note that the lower bound on $m_t$ \eqref{eq:nonconvex-m-lower} in the analysis of the nonconvex theorem still holds. Therefore, we have
		\begin{align*}
			G_t^{-\frac74}\delta_t &\overset{(i)}{\geq}\frac{\theta^2p_{m_t}G_t^{\frac14}}{8L} \overset{(ii)}{\geq}\frac{\theta^2m_t^2G_t^{\frac14}}{5L} \overset{(iii)}{\geq}\frac{\theta^2m_tG_t^{\frac14}}{10L}
			\frac{\bigl(L^2/MG_t\bigr)^{\frac14}} {\sqrt{\log_2 (\frac{4L^2}{MG_t})}} \overset{(iv)}{\geq} \frac{Cm_t}{L^{\frac12}M^{\frac14}\sqrt{\log_2 (1\!+\!\frac{L^2}{M\epsilon})}}. 
		\end{align*}
		Here $(i)$ is by Proposition~\ref{prop:convex-block}, $(ii)$ is by $p_{m_t}\geq1.6m_t^2$, $(iii)$ is by \eqref{eq:nonconvex-m-lower}, and $(iv)$ is by $G_t>\epsilon$.  
		
		On the other hand, the convexity and $L$-smoothness of $f$ implies $\epsilon<G_t\leq LR$ and $\Delta_t\leq RG_t$, which further gives $G_t\geq\max\{\epsilon,\Delta_t/R\}$. Therefore, we have 
		\begin{align*}
			\sum_{t=0}^{K-1}&G_t^{-7/4}\delta_t \leq \sum_{t=0}^{K-1}\frac{\Delta_t-\Delta_{t+1}}{\max \{\epsilon, \Delta_t/R\}^{7/4}} \leq \sum_{t=0}^{K-1} \int_{\Delta_{t+1}}^{\Delta_t}\frac{\dd s}{\max\{\epsilon,s/R\}^{7/4}}\\
			& \qquad =\int_{\Delta_K}^{\Delta_0}\frac{\dd s}{\max\{\epsilon,s/R\}^{7/4}} \leq \int_0^{R\epsilon}\!\!\! \epsilon^{-7/4}\dd s + R^{7/4}\!\! \int_{R\epsilon}^{\infty}\!\!s^{-7/4}\dd s = \frac73R\epsilon^{-3/4}.
		\end{align*}
		Combining the three main inequalities of this part gives
		\begin{equation}\label{eq:convex-paid-total}
			|\mathcal K_1|+\sum_{t\in\mathcal K_2}m_t \leq O\left( L^{1/2}M^{1/4}R\epsilon^{-3/4}\sqrt{\log_2 \!\left(1+\frac{L^2}{M\epsilon}\right)}\right).
		\end{equation} 
		
		\noindent\textbf{Epochs in $\mathcal K_3$ followed by $\mathcal{K}_2$.}\,\,
		As shown in the proof of Theorem~\ref{thm:nonconvex-main}, a consecutive sequence of epochs in $\mathcal K_3$ cannot be followed by an epoch in $\mathcal K_1$, and requires at most $C m_t$ updates if followed by $t\in\mathcal K_2$. Thus \eqref{eq:convex-paid-total} also bounds their total number of updates. It remains to bound epoch $K$ and the final consecutive sequence of epochs in $\mathcal K_3$ preceding it.\vspace{0.2cm}
		
		\noindent\textbf{Epochs in $\mathcal K_3$ followed by (and include) final epoch $K$.}\,\,
		Suppose $k\geq1$ iterations have been executed in epoch $t$. The first $k-1$ steps are all unclipped steps as the clipping ends an epoch. For $k\geq2$, the decomposition of $[a_{t,0},\ldots,a_{t,k-2}]$ in Lemma~\ref{lem:preceding-structure}(i) contains at most one sub-block of each length, hence one of length is at least $(k-1)/2$. Summing the descent inequality in Proposition~\ref{prop:convex-descent} over these iterations and applying \eqref{eq:general-mass} gives
		$$
		\Delta_t \geq f(x_{t,0})-f(x_{t,k-1}) \geq \frac{1}{4L}\sum_{i=0}^{k-2}a_{t,i}\|g_{t,i}\|^2 \geq \frac{\theta^2G_t^2}{4L}\sum_{i=0}^{k-2}a_{t,i} \geq \frac{\theta^2G_t^2(k-1)^2}{16L}.
		$$
		where the last inequality is by \eqref{eq:general-mass} and the above decomposition result. For $k=1$, the inequality holds trivially.  Again, using $G_t\leq LR$ and $\Delta_t\leq RG_t$, we obtain
		$$ k \leq 1 + \frac{4\sqrt{L\Delta_t}}{\theta G_t} \leq \left(1+\frac4\theta\right)\sqrt{\frac{LR}{G_t}}.
		$$
		Each epoch in the final consecutive sequence in $\mathcal K_3$ satisfies $G_{t+1}\leq\theta G_t$, while every epoch starts with $G_t>\epsilon$. Summing this bound backwards over these epochs and epoch $K$ gives 
		$$ \left(1+\frac4\theta\right)\sqrt{\frac{LR}{\epsilon}} \sum_{j=0}^{\infty} \theta^{j/2} \leq O\left(\sqrt{\frac{LR}{\epsilon}}\right).
		$$
		Combining this bound with \eqref{eq:convex-paid-total} proves this theorem.
	\end{proof}
	
\section{Numerical Experiments}
\label{sec:numerical}

In this section, we assess the numerical performance of NACC on squared-loss matrix completion, one-bit matrix completion, and neural network. The experiments were run on a MacBook Pro equipped with an Apple M5 Pro chip and $24$ GB of memory.

\subsection{Experimental Setup}
\label{sec:experimental-setup}

We compare NACC with four accelerated first-order methods, using the following implementations.

\begin{itemize}
    \item \textbf{Ada-RAGD-NC} \citep{li2023restarted} is an adaptive restarted accelerated gradient method. We follow the practical implementation in Section~5 of \citet{li2023restarted}\footnote{\url{https://github.com/lihuanML/RestartAGD}} with the internal tolerance $\epsilon=10^{-4}$. The initial stepsizes for the individual problems are given in the later discussions.

    \item \textbf{MT-RAGD} \citep{marumo2024parameterfree} is a parameter-free restarted accelerated gradient method. As suggested by \citet{marumo2024parameterfree}, we set $(L_{\mathrm{init}},M_0,\alpha,\beta) =(10^{-3},10^{-16},2,0.9)$.

    \item \textbf{MT-RHB} \citep{marumo2025universal} is a universal restarted heavy-ball method that does not require an initial Hessian Lipschitz estimate. Following \citet{marumo2025universal}, we set $(\ell_{\mathrm{init}},\alpha,\beta)=(10^{-3},2,0.1)$. The implementations of MT-RAGD and MT-RHB are available online.\footnote{\url{https://github.com/n-marumo/restarted-hb}}

    \item \textbf{CUPG} \citep{carmon2017convex} is an accelerated gradient method with explicit negative-curvature detection and exploitation. We follow the numerical implementation used by \citet{li2023restarted}.

    \item \textbf{NACC} is our proposed method. In the practical implementation, we increase the energy budget $E_t^{\max}$, restrict the schedule length $m$, and adaptively adjust the estimates $\widehat L$ and $\widehat M$ of the gradient and Hessian Lipschitz constants $L$ and $M$, respectively. Across all experiments, we fix $\theta=0.9$.
\end{itemize}

All methods start from the same initial point for each problem. For each experiment, we compare the behavior of objective values and gradient norms over a fixed number of iterations and a fixed amount of wall-clock time.

\subsection{Matrix Completion with Squared Loss}
\label{sec:matrix-completion-square}

We first consider squared-loss matrix completion \citep{negahban2012restricted,hardt2014understanding} on MovieLens 10M and MovieLens 20M \citep{harper2015movielens}, and Netflix \citep{bennett2007netflix}, whose  dimensions are given in Table~\ref{tab:rating-data}.  Let $A\in\mathbb R^{n_1\times n_2}$ denote the rating matrix, $\Omega$ the set of observed entries, and $N=|\Omega|$. The problem is formulated in Section~5.1 of \citet{li2023restarted} as
\begin{equation}
\label{eq:experiment-square}
    \min_{U\in\mathbb R^{n_1\times r},\,
          V\in\mathbb R^{n_2\times r}}
    f(U,V)
    =\frac{1}{2N}\sum_{(i,j)\in\Omega}
        \bigl(u_i^\top v_j-A_{ij}\bigr)^2
     +\frac{1}{2N}\|U^\top U-V^\top V\|_F^2,
\end{equation}
where $u_i^\top$ and $v_j^\top$ are the rows of $U$ and $V$. We use all observed ratings, set $r=10$, and initialize every method from the same rank-$10$ sparse SVD.

\begin{table}[t]
    \centering
    \caption{Rating datasets used in matrix completion experiments. The number of optimization variables is $d=r(n_1+n_2)$ with $r=10$.}
    \label{tab:rating-data}
    \small
    \begin{tabular}{lrrrr}
        \toprule
        Dataset & $n_1$ & $n_2$ & $N$ & $d$ \\
        \midrule
        MovieLens 10M & 10,677 & 69,878 & 10,000,054 & 805,550 \\
        MovieLens 20M & 26,744 & 138,493 & 20,000,263 & 1,652,370 \\
        Netflix & 17,770 & 480,189 & 100,480,507 & 4,979,590 \\
        \bottomrule
    \end{tabular}
\end{table}

For Ada-RAGD-NC, we use the dataset-specific initial stepsizes from the public code of \citet{li2023restarted}. In NACC, we  take $\widehat L=1/\eta$, and $\widehat M=10^{-5}\|g_2-2g_1+g_0\|/\delta^2$, where
\begin{equation}
\label{eq:experiment-probe}
    \begin{aligned}
        g_0=\nabla f(x_0),\quad
        s=-\delta\frac{g_0}{\|g_0\|},\qquad \delta=10^{-3},\quad
        g_1=\nabla f(x_0+s),\quad
        g_2=\nabla f(x_0+2s).
    \end{aligned}
\end{equation}

\begin{figure*}[h]
    \centering
    \includegraphics[width=\textwidth]{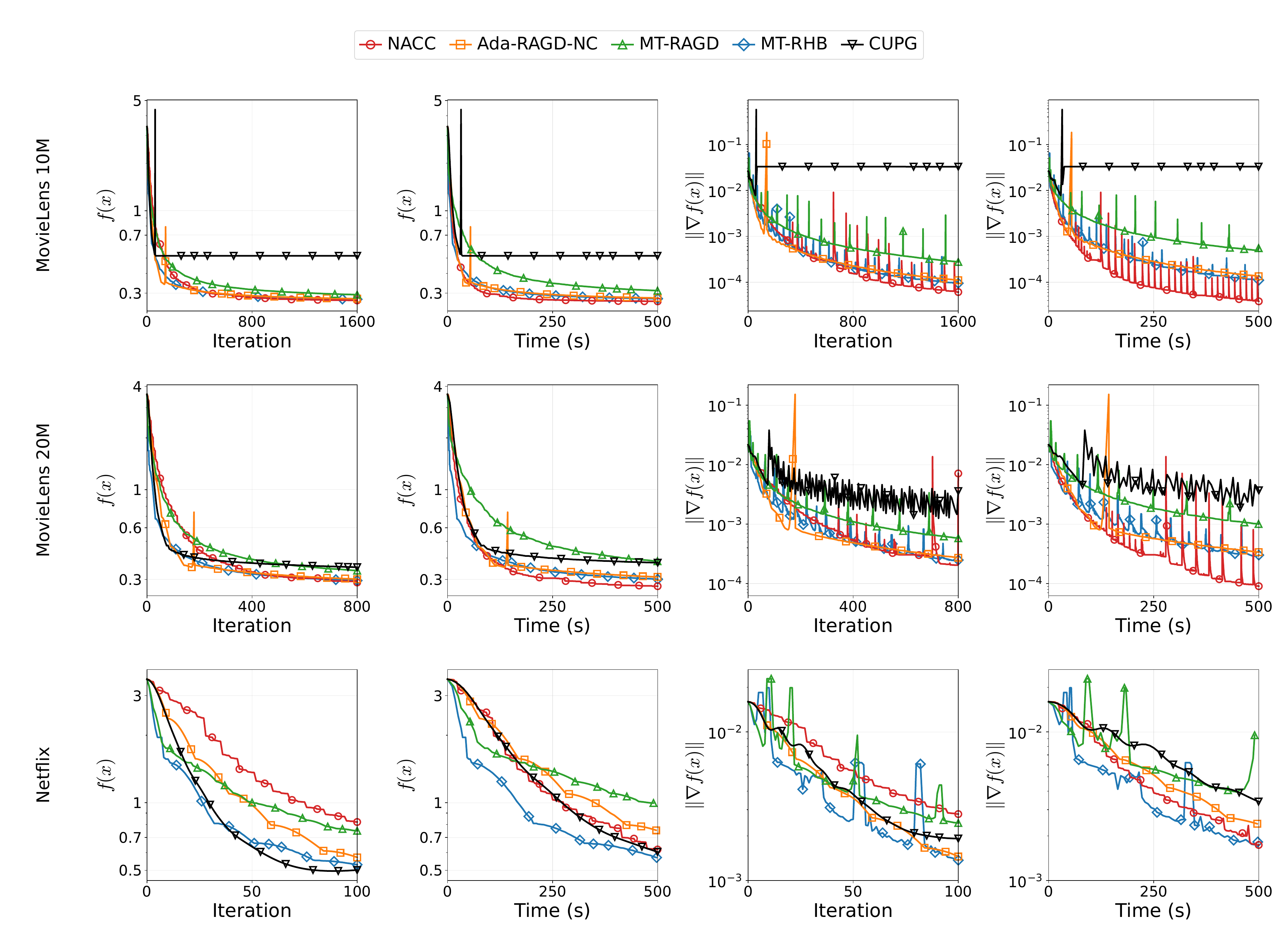}
    \caption{Squared-loss matrix completion on MovieLens 10M (top row), MovieLens 20M (middle row), and Netflix (bottom row). From left to right, the columns show the objective value versus iterations and time, followed by the gradient norm versus iterations and time.}
    \label{fig:matrix-square}
\end{figure*}

Figure~\ref{fig:matrix-square} reports the objective values and gradient norms of the compared methods against iteration count and elapsed time. Although NACC is not consistently ahead in the iteration plots, its advantage is more pronounced when progress is measured by elapsed time. In particular, NACC exhibits faster decay of the gradient norm over time on the two MovieLens datasets. It also achieves the lowest gradient norm on all three datasets at the end of the $500$-second comparison window. This advantage may be attributed to its simple gradient-descent updates, which allow it to complete more iterations within the same amount of time.



\subsection{One-Bit Matrix Completion}
\label{sec:matrix-completion-onebit}

We next consider one-bit matrix completion \citep{davenport2014onebit}, using a logistic loss on the same three rating datasets. Writing $\bar A=N^{-1}\sum_{(i,j)\in\Omega}A_{ij}$ for the mean observed rating, we set $b_{ij}=1$ when $A_{ij}>\bar A$ and $b_{ij}=-1$ otherwise. Following \citet{li2023restarted}, we minimize
\begin{equation}
\label{eq:experiment-onebit}
    \min_{U\in\mathbb R^{n_1\times r},\,
          V\in\mathbb R^{n_2\times r}}
    f(U,V)
    =\frac{1}{N}\sum_{(i,j)\in\Omega}
        \log\bigl(1+\exp(-b_{ij}u_i^\top v_j)\bigr)
     +\frac{1}{2N}\|U^\top U-V^\top V\|_F^2.
\end{equation}
As in the previous experiment, we set $r=10$ and use the same spectral initialization procedure, applied here to the sparse matrix whose $(i,j)$ entry is $b_{ij}$ if $(i,j)\in\Omega$ and zero otherwise. We use the  stepsizes $\eta_0$ from the public Ada-RAGD-NC implementation and initialize NACC with $\widehat L=1/\eta_0$ and $\widehat M$ from \eqref{eq:experiment-probe}, evaluated for \eqref{eq:experiment-onebit}.
\begin{figure*}[h]
    \centering
    \includegraphics[width=\textwidth]{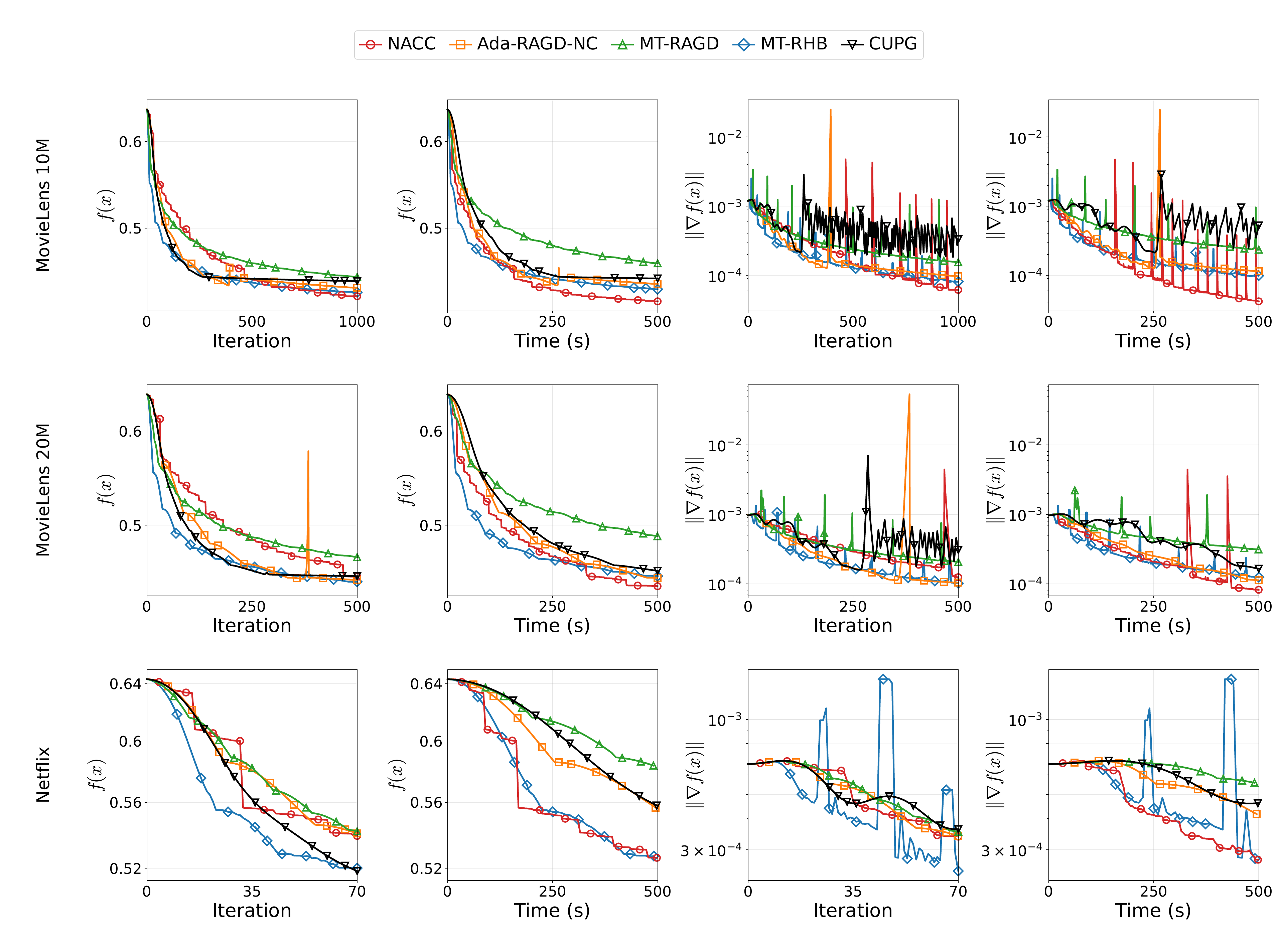}
    \caption{One-bit matrix completion on MovieLens 10M (top row), MovieLens 20M (middle row), and Netflix (bottom row). From left to right, the columns show the objective value versus iterations and time, followed by the gradient norm  versus iterations and time.}
    \label{fig:matrix-onebit}
\end{figure*}

Figure~\ref{fig:matrix-onebit} reports the objective values and gradient norms for one-bit matrix completion against iteration count and elapsed time. As in the previous experiment in Section~\ref{sec:matrix-completion-square}, NACC's advantage is more pronounced when progress is measured by elapsed time, and achieves the lowest objective value and gradient norm on all three datasets at the end of the $500$-second comparison window.\\

In both the squared-loss and one-bit matrix completion experiments, several sharp increases in the gradient norm coincide with substantial decreases in the objective value on MovieLens 10M and 20M. One possible explanation is overshooting along positive-curvature directions caused by large Chebyshev steps, which may temporarily amplify gradient components even as the objective decreases. Negative curvature may also contribute to this phenomenon, as it can amplify gradient components while inducing objective decrease. Our analysis accounts for this effect by controlling the negative-curvature components through a backward estimate (Lemma~\ref{lem:negative-gradient-bound}). Unlike CUPG \citep{carmon2017convex}, which explicitly detects and exploits negative curvature, NACC can benefit from it through ordinary gradient updates without a separate detection procedure.



\subsection{Neural Network Classification}
\label{sec:neural-network-classification}
We finally turn to full-batch classification on MNIST \citep{lecun1998gradient}. We adopt the network architecture and cross-entropy objective of \citet{marumo2025universal}, replacing sigmoid with ReLU in both hidden layers. The fully connected network has widths $784$--$32$--$16$--$10$, includes biases in all affine layers, and contains $d=25,818$  parameters. We select $N=10,000$ training images by class-stratified sampling.

For an image $a$ and parameter vector $x$, let $z(a;x)\in\mathbb R^{10}$ denote the network logits, and let $y_i\in\{0,\ldots,9\}$ be the label of image $a_i$. We consider
\begin{equation}
\label{eq:experiment-classification}
    \min_{x\in \mathbb R^d} f(x)=\frac{1}{N}\sum_{i=1}^N
    \left[
        \log\left(\sum_{c=0}^{9}\exp(z_c(a_i;x))\right)
        -z_{y_i}(a_i;x)
    \right].
\end{equation}
We apply \eqref{eq:experiment-probe} to initialize \( \widehat L={\|g_1-g_0\|}/{4\delta}, \; \widehat M=10^{-5}{\|g_2-2g_1+g_0\|}/{\delta^2} \) in NACC. For Ada-RAGD-NC, the initial stepsize is $\eta_0=1/\widehat L$.

\begin{figure*}[t]
    \centering
    \includegraphics[width=\textwidth]{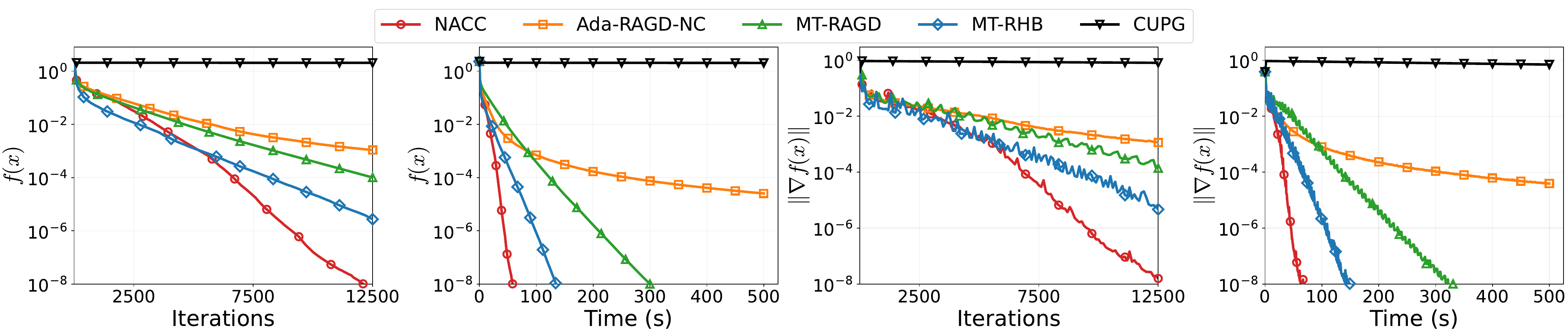}
    \caption{MNIST classification with a ReLU network. Objective values (left two panels) and gradient norms (right two panels) versus iterations and elapsed time.}
    \label{fig:mnist-classification}
\end{figure*}

Figure~\ref{fig:mnist-classification} shows that NACC significantly outperforms all competing methods in reducing both the gradient norm and objective value, in terms of both iteration count and wall-clock time. In particular, NACC reaches a gradient norm of $10^{-8}$ in $63.1$ seconds, compared with $141.6$ seconds for MT-RHB and $319.0$ seconds for MT-RAGD, while Ada-RAGD-NC and CUPG fail to reach this threshold within $500$ seconds. These results highlight the effectiveness of NACC and its potential for broader applications in neural network optimization.

Overall, the numerical results demonstrate the practical effectiveness of NACC across different nonconvex optimization problems. For matrix completion, NACC achieves substantially better wall-clock performance, particularly on the MovieLens datasets, although it does not always outperform the competing methods in terms of iteration count. For neural network classification, NACC exhibits clear advantages in reducing both the objective value and gradient norm, measured against either iterations or elapsed time. These results highlight the potential of adaptive Chebyshev stepsize scheduling for practical nonconvex optimization, showing that significant acceleration can be achieved while retaining the simple gradient-descent update without momentum or extrapolation.

\appendix

	\section{Chebyshev Schedule and Polynomial Estimates}
	\label{sec:cheb-toolkit}
	This section only contains the technical supporting lemmas of the Chebyshev step schedule itself, and is not directly related to the optimization analysis. We isolate its presentation out for better readability for optimization readers. 
	
	\paragraph{Basics of Chebyshev schedule.}
	For any $m\geq 0$, let $\cT_m$ be the first kind Chebyshev polynomial, and let $\cU_m$ be the second kind Chebyshev polynomial. We use the following standard identities
	\begin{equation}\label{eq:cheb-standard}
		\cT_m(\cos v)=\cos(mv),
		\qquad
		\cU_m(\cos v)=\frac{\sin((m+1)v)}{\sin v},
		\qquad
		\cT_{2m}(z)=2\cT_m(z)^2-1.
	\end{equation}
	Recall that in Section \ref{sec:prelim}, for any ordered sequence $Q=[a_1,\ldots,a_{m}]$, we define its associated polynomial as $P_Q(z) = \prod_{i=1}^{m} (1-a_i z)$. When $Q = \Acal_{m,\phi}$ is the Chebyshev block of size $m$, denote $p_Q =p_{m,\phi}$ the largest element in $\Acal_{m,\phi}$, we define the associated \emph{preparation polynomial} by
	\begin{equation}\label{eq:PQ-FQ-def}
		F_Q(z):=\frac{P_Q(z)}{1-p_Qz}.
	\end{equation}
	For the Chebyshev block that adopts the specific reversed fractal ordering in this paper, $p_Q$ is always the last element of the block. Therefore, for $Q = \Acal_{m,\phi} = [a_1,\ldots,a_{m}]$, we can also write $F_Q(z) = \prod_{i=1}^{m-1} (1-a_i z)$, equivalently. Based on these definition, we present a few handy properties of the reversed fractal schedule that will be used in later discussion. 
	
	First, we state two structural properties that follow from the recursive construction of the fractal schedule \citep{lebedev1971ordering}, which was systematically summarized and presented by \citet{agarwal2021acceleration}. In this paper, we slightly modify and adapt their analysis to the reversed ordering in \eqref{eq:A-def}. A short proof is included here for completeness. 
	
	\begin{lemma}\label{lem:preceding-structure}
		Let $\Acal_{m,\phi} = [a_1,a_2,\cdots,a_m]$ for $m\in\{2^k:k\in\mathbb{N}\}$ and $\phi\in[\pi/2,\pi)$, we have\vspace{0.1cm} \\ 
		\noindent\textbf{(i). } For any $2\leq r\leq m$, we can partition $[a_1,\cdots,a_{r-1}]$ into $k\geq1$ sub-blocks, each block has the form $\Acal_{m',\phi'}$ with $m'\in\{2^k:k\in\mathbb{N}\}$ and $\phi'\in[\pi/2,\pi)$.\vspace{0.1cm} \\ 
		\noindent\textbf{(ii). } For any $1\leq r\leq m$, the largest recursive sub-block $\Acal_{m',\phi'}$ with $m'\in\{2^k:k\in\mathbb{N}\}$ ending at $a_r$ will satisfy $\phi'\in[\pi/2,\pi)$, and $a_r\leq 2(m')^2$. 
	\end{lemma}
	\begin{proof}
		For (i), the partition can be achieved by the following procedure. We start from $a_{r}$ and travel upwards in the binary tree that recursively defines the sequence, see Fig \ref{fig:tree-m4} for example. Each time, if we arrive at a parent node from a left child, we keep move upward. If we arrive at a parent node from a right child, we collect the leaves of the sub-tree rooted at the left child of this parent node as one sub-block, and we repeat this until we reach the root node $\Acal_{m,\phi}$. As these sub-blocks corresponds to some left child node, according to the definition \eqref{eq:A-def}, each block $\Acal_{m',\phi'}$ will naturally have $\phi'=\pi-\phi''/2\in[\pi/2,\pi)$ for some $\phi''\leq \pi$, and $m'\in\{2^k:k\in\mathbb{N}\}$ as it contains the leaves of a binary sub-tree. 
		
		For (ii), as mentioned in the proof of Lemma~\ref{lem:curvature-weighted-gradient}, the largest recursive sub-block with $a_r$ being the right most leaf can be constructed as follows. If $a_r$ is already a left child, then the block is $\{a_r\}$ alone. Otherwise, we search the binary tree upwards until we meet a parent that is a left-child of the the parent (or if it is already the root when $r=m$). According to the proof of (i), we know $\phi'\in[\pi/2,\pi)$ for this sub-block. Then $a_r = p_{m',\phi'}\leq p_{m',\pi/2} \leq 2(m')^2$ holds as discussed in Section~\ref{sec:prelim}.  
	\end{proof}
	
	Next, we establish a few new properties that are not adapted from \citet{lebedev1971ordering,agarwal2021acceleration}. They are also important tools in our analysis.  
	
	\begin{lemma}
		\label{lem:tail-preparation}
		Let $Q=\Acal_{m,\phi}=[a_1,\ldots,a_m]$ with $m\in\{2^k:k\in\mathbb N\}$ and $\phi\in(0,\pi)$. Then  
		$$
		\abs{P_{Q_{i:m-1}}(z)}\leq1, \quad \forall z\in[0,1],\quad \forall 1\leq i\leq m.
		$$
		If $\phi\in[\pi/2,\pi)$, for  $\forall m'\in\{2^k:k\in\mathbb N\}, m'<m$, the sub-block $Q':=Q_{m-m'+1:m}$ satisfies
		$$
		\abs{F_{Q'}(z)}\leq1, \qquad \abs{P_{Q'}(z)}\leq2\left({m}/{m'}\right)^2, \qquad  \forall z\in[0,1].
		$$
	\end{lemma}
	
	Before presenting the proof, it is worth noting that for the second inequality, we actually have $Q' = \Acal_{m', \phi'}$ with $\phi'=\frac{m'}{m}\phi\in(0,\pi/2)$, hence $|F_{Q'}(z)|\leq 1, \forall z\in[0,1]$ cannot be implied by the known result Lemma \ref{lemma:prefix-stability}, which requires $\phi'\in[\pi/2,\pi)$. Instead it is a consequence of the first inequality of the lemma, by setting $i=m-m'+1$, we specifically single it out for the convenience of later reference. \vspace{0.2cm}
	
	\begin{proof}
		First, let us prove the bound for $|P_{Q_{i:m-1}}|$. For this inequality, let us prove it by induction. If $m=1$, then $Q_{i:m-1} = \varnothing$ and $P_\varnothing(z)=1$ by default. For $m\geq2$, definition \eqref{eq:A-def} gives $Q=\big(\Acal_{m/2, \pi-\phi/2}, \Acal_{m/2,\phi/2}\big).$ Assume the inequality already holds for $m/2$, we prove it for $m$.  If $i>m/2$, then $Q_{i:m-1}$ is also a sub-block of $\Acal_{m/2,\phi/2}$ and hence $|P_{Q_{i:m-1}}|\leq 1$ for all $z\in[0,1]$ according to our induction hypothesis for $m/2$. If $i\leq m/2$, then the induction hypothesis implies
		$$
		\begin{aligned}
			\abs{P_{Q_{i:m-1}}(z)} = \abs{P_{Q_{i:m/2-1}}(z)} \cdot \abs{(1-a_{m/2} z)F_{\Acal_{m/2,\phi/2}}(z)} \leq \abs{(1-a_{m/2} z)F_{\Acal_{m/2,\phi/2}}(z)}.
		\end{aligned}
		$$
		Set $r_1=\sin^2(\frac{\phi}{2m})$ and $r_2=\sin^2(\frac{2\pi-\phi}{2m}))=\frac{1}{a_{m/2}}$. Then for $0\leq z\leq r_2$, we have 
		$$\big|(1-a_{m/2}z)F_{\Acal_{m/2,\phi/2}}(z)\big| \leq \big|(1-a_{m/2}z)\big|\cdot \big|F_{\Acal_{m/2,\phi/2}}(z)\big|\leq 1$$
		because both terms in the product lie in $[0,1]$. For $r_2< z\leq 1$,  we have 
		$$
		\abs{(1-a_{m/2}z)F_{\Acal_{m/2,\phi/2}}(z)} \leq\frac{r_1}{r_2} \cdot \cot^2\Big(\frac{\phi}{4}\Big) \cdot \frac{z-r_2}{z-r_1} \leq 1,
		$$
		where the first inequality is due to \eqref{eq:PQ-cheb}, and the second inequality is by $r_1<r_2$ and the fact that $\sin(\frac{\phi}{2m})/\sin(\frac{2\pi-\phi}{2m})\leq\tan(\frac{\phi}{4})$ for $\phi\in(0,\pi)$.  Conbining the above arguments finishes the induction.
		
		Next, we prove the inequalities for the sub-block $Q'$. According to the  recursive construction definition \eqref{eq:A-def}, we know $Q'=\Acal_{m',\phi'}$ and $\phi':=\frac{m'}m\phi\in(0,\pi/2)$. As commented before the proof, the bound on $F_{Q'}$ is directly implied by the first inequality. We only need to prove the inequality for $P_{Q'}$, for which  \eqref{eq:PQ-cheb} gives
		$$
		\abs{P_{Q'}(z)} = \bigg|\frac{\cT_{m'}(1-2z)-\cos\phi'}{1-\cos\phi'}\bigg|
		\leq\frac{1+\cos\phi'}{1-\cos\phi'}
		=\cot^2(\phi'/2)
		\leq\frac{4}{(\phi')^2}
		=\frac{4}{(\frac{m'}{m}\phi)^2}<2\left(\frac{m}{m'}\right)^2,
		$$
		where the last inequality is by $\phi\in[\pi/2,\pi)$.
	\end{proof}
	
	Building on the preceding preparation bound, we establish a few handy   inequalities for the preparation polynomial $F_Q(z)$, see definition in \eqref{eq:PQ-FQ-def}. These inequalities will be used in the analysis of nonconvex and convex cases, respectively.  
	
	\begin{lemma}\label{lem:preparation-moment}
		Let $Q=\Acal_{m,\phi}=[a_1,\ldots,a_{m}]$, $\phi\in[\pi/2,\pi)$, then we have $\sqrt{a_mz}\,\abs{F_Q(z)}\leq\frac76$ for all $z\in\big[0,\frac12\big]$. 
	\end{lemma}
	\begin{proof}
		For $m=1,2$, this inequality can be verified by direct computation, and is hence omitted. For $m\geq4$, when $a_mz\leq 1$, Lemma~\ref{lem:tail-preparation} immediately implies $\sqrt{a_mz}\abs{F_Q(z)}\leq\sqrt{a_mz}\leq1$. 
		If $\sqrt{a_mz}\geq38/25$, Lemma~\ref{lemma:prefix-stability} gives
		\begin{equation}
			\sqrt{a_mz}\abs{F_Q(z)}
			=\frac{\sqrt{a_mz}\abs{P_Q(z)}}{a_mz-1}
			\leq\frac{\sqrt{a_mz}}{a_mz-1}
			\leq \max_{u\geq\frac{38}{25}}\left\{\frac{u}{u^2-1}\right\}\leq \frac76.\nonumber
		\end{equation}
		Finally, it remains to verify the interval $1<\sqrt{a_mz}<38/25$. Write $z=1/\csc^2v = \sin^2v$. Then $a_m z = \csc^2(\phi/2m)z>1$ implies with $v>\frac{\phi}{2m}$. By \eqref{eq:PQ-cheb}, we have 
		$$
		\sqrt{a_mz}\abs{F_Q(z)}
		=\frac{\sin v\sin(\frac{\phi}{2m})\,
			\abs{\sin(mv+\frac{\phi}{2})\sin(mv-\frac{\phi}{2})}}
		{\sin^2(\frac{\phi}{2})\sin(v+\frac{\phi}{2m})\sin(v-\frac{\phi}{2m})}\leq\frac{m\sin v\sin(\frac{\phi}{2m})}
		{\sin^2(\frac{\phi}{2})\sin(v+\frac{\phi}{2m})},
		$$
		where we used $\abs{\sin(mv-\frac{\phi}{2})}\leq m\sin(v-\frac{\phi}{2m})$. Then by $m\sin\big(\frac{\phi}{2m}\big)\leq\frac{\phi}{2}\leq\frac{\pi}{2}\sin^2\big(\frac{\phi}{2}\big)$ for any $\phi\in[\pi/2,\pi)$, and  
		$\sin\big(v+\frac{\phi}{2m}\big)\geq\cos v \sin v+\cos v\sin\big(\frac{\phi}{2m}\big)$ for $v>\phi/(2m)$, we have    
		\begin{equation} 
			\sqrt{a_mz}\abs{F_Q(z)}
			\leq\frac\pi2\frac{\sin v}{\bigl(\sin v+\sin(\phi/(2m))\bigr)\cos v}
			=\frac\pi2\frac{\sqrt{a_mz}}{(\sqrt{a_mz}+1)\sqrt{1-z}}. \nonumber
		\end{equation}
		Finally, as $m\geq4$ and $\sqrt{a_mz}<38/25$, we know $z<(\frac{38}{25}\sin(\frac{\phi}{2m}))^2\leq (\frac{38}{25}\sin(\frac{\pi}{8}))^2:=z_0$, by numerical computation, we have 
		$$\sqrt{a_mz}\abs{F_Q(z)} \leq  \sup_{1<u<\frac{38}{25}}\left\{\frac{\pi u}{2\sqrt{1-z_0}(u+1)}\right\}\leq \frac{7}{6}.$$
		Combining the above arguments proves the lemma.  
	\end{proof}

	\begin{lemma}\label{lem:dissipative}
		Let $Q=\Acal_{m,\phi}=[a_1,\ldots,a_m]$ with $m\geq2$ and $\phi\in[\pi/2,\pi)$. Then
		\begin{equation}
			\left(1+\frac{a_mz}{6}\right)F_Q(z)^2\leq1,\qquad \forall z\in[0,1/2]. \nonumber
		\end{equation}
	\end{lemma}
	\begin{proof}
		If $a_mz\leq1$, then $0\leq1-a_jz\leq1$ for every $j$ as $a_j\leq a_m$. Therefore, we have 
		$0\leq F_Q(z)=\prod_{j=1}^{m-1}(1-a_jz)\leq1-a_{m/2}z,$ which further implies that  
		$$
		\begin{aligned}
			1\!-\!F_Q(z)^2 \geq 1\!-\!(1-a_{m/2}z)^2 \geq 2a_{m/2}z(1\!-\!a_{m/2}z)^2
			\geq 2a_{m/2}zF_Q(z)^2
			\overset{(i)}{\geq}\frac{a_mz}{6}F_Q(z)^2,
		\end{aligned}
		$$
		rearranging the terms proves the target inequality. In the above inequalities, we used the following fact for $\phi\in[\pi/2,\pi)$ and $m\geq2$:   
		$$
		\frac{a_{m/2}}{a_m}
		=\frac{\sin^2(\frac{\phi}{2m})}{\sin^2(\frac{2\pi-\phi}{2m})}
		\geq\left(\frac{\phi}{2\pi-\phi}\right)^2\geq\frac19.
		$$
		The second inequality above uses $1-(1-u)^2-2u(1-u)^2=u^2(3-2u)\geq0$ for $u\in[0,1]$.
		On the other hand, if $a_mz>1$, then Lemma~\ref{lem:preparation-moment} implies that $\frac{a_mz}{6}F_Q(z)^2\leq 49/216.$
		To bound the remaining $F_Q(z)^2$, let us again write $z = \sin^2 v$ and decompose $Q_{1:m-1}$ into the blocks $\Acal_{m',\pi-(m'/m)\phi}$ for $m'=m/2,m/4,\ldots,1$. Since $mv>\phi/2\geq\pi/4$ and $v\leq\pi/4$, one of the $m'$s will satisfy $\pi/8\leq m'v\leq\pi/4$. All these blocks have $\pi-(m'/m)\phi\in[\pi/2,\pi)$, and hence Lemma~\ref{lemma:prefix-stability} and \eqref{eq:PQ-cheb} imply
		$$ 
		\abs{F_Q(z)} \leq \abs{P_{\Acal_{m',\pi-(m'/m)\phi}}(z)} = \frac{\cos(2m'v)+\cos(\frac{m'\phi}{m})}{1+\cos(\frac{m'\phi}{m})} \overset{(i)}{\leq}\frac{1+1/\sqrt2}{2}, 
		$$
		where (i) id due to  $0\leq\cos(2m'v)\leq1/\sqrt2$. Hence, $\big(1+\frac{a_mz}{6}\big)F_Q(z)^2
		\leq\big(\frac{1+1/\sqrt2}{2}\big)^2+\frac{49}{216}\leq 1$.\vspace{0.3cm}
	\end{proof}

	\begin{lemma}\label{lem:small-phase}
		Let $Q=\Acal_{m,\psi}$ with $\psi\in(0,\pi/2]$, then $\sqrt z\,\abs{F_Q(z)}\leq\frac{\pi}{2m}$ for all $z\in[0,1]$. 
	\end{lemma}
	\begin{proof}
		By Lemma~\ref{lem:tail-preparation}, we know $\abs{F_Q(z)}\leq1$. If   $z\leq2\sin^2(\psi/2m)$, then  
		$$\sqrt{z}\abs{F_Q(z)} \leq \sqrt{2}\sin(\psi/2m)\leq\pi/2m.
		$$
		Otherwise, for for $z\geq2\sin^2(\psi/2m)$, equation \eqref{eq:PQ-cheb} gives  $$\abs{F_Q(z)}=\frac{\sin^2(\psi/2m)\,\abs{\cT_m(1-2z)-\cos\psi}}{\bigl(z-\sin^2(\psi/2m)\bigr)(1-\cos\psi)} \overset{(i)}{\leq} \frac{2\sin^2(\psi/2m)} {z\sin^2(\psi/2)} \overset{(ii)}{\leq} \frac{\pi^2}{4m^2z}.$$
		Here (i) is by $\abs{\cT_m(1-2z)-\cos\psi}\leq 2$, $z-\sin^2(\psi/2m)\geq z/2$, and $1-\cos\psi = 2\sin^2(\psi/2)$. And (ii) is due to $\sin(\psi/2m)\leq \psi/2m$ and $\sin(\psi/2)\geq \sqrt{2}\psi/\pi$.  Consequently, we obtain $zF_Q(z)^2\leq z\abs{F_Q(z)} \leq\frac{\pi^2}{4m^2}$ as desired. \vspace{0.3cm} 
	\end{proof}

	\begin{lemma}\label{lem:weighted-response}
		Let $Q=\Acal_{m,\phi}=[a_1,\ldots,a_m]$ with $\phi\in[\pi/2,\pi)$. For any $0\preceq B\preceq I/2$, we have 
		\begin{equation}
			\sum_{j=1}^{m-1}a_j
			\left\|\sqrt{a_mB}\,P_{Q_{j+1:m-1}}(B)\right\|
			\leq\frac{\pi}{\sqrt2}m^2\log_2m. \nonumber
		\end{equation}
	\end{lemma}
	\begin{proof}
		For each $j<m$, let $m'$ be the size of the largest recursive sub-block ending at $a_j$, denoted by $\Acal_{m',\phi'}, \phi'\in[\pi/2,\pi)$, its construction can be found in the proof of Lemma \ref{lem:preceding-structure}.  Its right sibling $Q_{j+1:j+m'} = \Acal_{m',\phi''}$ satisfies $\phi''=\pi-\phi'\in(0,\pi/2)$. As a result, for $z\in[0,1/2]$, we have 
		$$
		\sqrt z\,\abs{P_{Q_{j+1:m-1}}(z)}
		=\sqrt z\,\abs{F_{Q_{j+1:j+m'}}(z)}\,
		\abs{P_{Q_{j+m':m-1}}(z)}
		\leq\frac{\pi}{2m'},
		$$
		where the inequality is due to Lemma~\ref{lem:small-phase} and Lemma~\ref{lem:tail-preparation}. Applying the above inequality to the eigenvalues of $B$ gives $\|\sqrt{a_mB}\,P_{Q_{j+1:m-1}}(B)\| \leq \frac{\sqrt{a_m}\pi}{2m'}$. For each $m'\in\{2^k:k\in\mathbb{N}\}, m'<m$, there are $m/2m'$ such indices $j$ for which the largest recursive sub-block ending at $a_j$ has the same size $m'$, and for all these $j$'s,  $a_j\leq2(m')^2$ by Lemma~\ref{lem:preceding-structure}(ii). Their total contribution in the summation is therefore at most
		$$
		\frac{m}{2m'}\cdot2(m')^2\cdot\sqrt{a_m}\frac{\pi}{2m'}
		=\frac{\pi m\sqrt{a_m}}{2}
		\leq\frac{\pi}{\sqrt2}m^2.
		$$
		Summing over the $\log_2m$ possible $m'$s proves the lemma.
	\end{proof}

	\begin{lemma}\label{lem:ordered-sensitivity}
		Let $Q=\Acal_{m,\phi}$ with $\phi\in[\pi/2,\pi)$. For symmetric matrices $0\preceq B,\widetilde B\preceq I$, we have
		\begin{equation}
			\norm{P_Q(B)-P_Q(\widetilde B)}
			\leq\frac{2m^2}{1-\cos\phi}\norm{B-\widetilde B}.\nonumber
		\end{equation}
	\end{lemma}
	\begin{proof}
		Let $X$ and $Y$ be symmetric matrices satisfying $-I\preceq X,Y\preceq I$. Using the three-term recurrence of $\cT_m$ and $\cU_m$, one can prove by induction the following identity
		\begin{equation}
			\label{eq:ordered-diff}
			\cT_m(X)-\cT_m(Y)
			=(X-Y)\cU_{m-1}(Y)
			+2\sum_{k=1}^{m-1}\cT_k(X)(X-Y)\cU_{m-1-k}(Y).
		\end{equation}
		As $-I\preceq X,Y\preceq I$, we have $\norm{\cT_k(X)}\leq1$ and $\norm{\cU_k(Y)}\leq k+1$.  Hence
		\begin{equation}
			\norm{\cT_m(X)-\cT_m(Y)}
			\leq\left(m+2\sum_{k=1}^{m-1}(m-k)\right)\norm{X-Y}
			=m^2\norm{X-Y}.\nonumber
		\end{equation}
		Applying this bound to $X=I-2B$ and $Y=I-2\widetilde B$ proves the lemma.
	\end{proof}
	
	\section{Technical Proof of Section \ref{sec:nonconvex}}
	\label{sec:nonconvex-technical}
	This appendix contains the technical proof of Section \ref{sec:nonconvex} that were relegated to appendix, which is essentially the inequality \eqref{eq:positive-prod-res}. Next, we first establish a few preliminary tools for proving this inequality and then prove this inequality itself.  
	
	\subsection{A Few Preliminary Tools for \eqref{eq:positive-prod-res}} 
	\begin{lemma}\label{lem:cheb-matrix-bounds}
		Let $B$ be a matrix (not necessarily symmetric) s.t. $\norm{B-B_0}\leq d$ for some $0\preceq B_0\preceq I/2$.  Define $\mathcal Q_{n,\phi}(B):=\cT_n(I-2B)-\cos\phi\,I$. Then for every $n\in\mathbb N$, we have
		\begin{equation}\label{eq:T-growth}
			\norm{\cT_n(I-2B)}\leq\ee^{2n^2d}.
		\end{equation}
		For any $\phi\in(0,\pi)$ and any $D$ satisfying $\norm{D-D_0}\leq d$ for some $0\preceq D_0\preceq I/2$, we have 
		\begin{equation}\label{eq:Q-diff}
			\norm{\mathcal Q_{n,\phi}(B)-\mathcal Q_{n,\phi}(D)}
			\leq2n^2\ee^{2n^2d}\norm{B-D}.
		\end{equation}
		If in addition that $\phi\in[\pi/2,\pi)$, we also have
		\begin{equation}\label{eq:frozen-good-new}
			\left\|\frac{\mathcal Q_{n,\phi}(B)}{1-\cos\phi}\right\|
			\leq\ee^{4n^2d}.
		\end{equation}
	\end{lemma}
	\begin{proof}
		All three inequalities are immediate for $n=0$. To prove \eqref{eq:T-growth} for $n\geq1$, set $X=I-2B$ and $X_0=I-2B_0$, then $0\preceq X_0\preceq I$ and $\norm{X-X_0}\leq2d$. By \eqref{eq:ordered-diff}, we have 
		$$
		\begin{aligned}
			\norm{\cT_n(X)} & = \Big\|\cT_n(X_0) + (X-X_0)U_{n-1}(X_0) + 2\sum_{k=1}^{n-1}\cT_k(X)(X-X_0) U_{n-1-k}(X_0)\Big\| \\
			& \leq 1+2nd+4d\sum_{k=1}^{n-1}(n-k)\norm{\cT_k(X)},
		\end{aligned}
		$$ 
		where the inequality is due to  \eqref{eq:cheb-standard}, which implies  $\norm{\cT_k(X_0)}\leq1$ and $\norm{U_k(X_0)}\leq k+1$.  
		Then applying the same identity to the scalars $1+2d$ and $1$ gives
		$$
		\cT_n(1+2d)=1+2nd+4d\sum_{k=1}^{n-1}(n-k)\cT_k(1+2d).
		$$
		Since $d\geq0$, comparing the above two inequalities and using induction on $n$ from $\cT_0=1$ gives $\norm{\cT_n(X)}\leq \cT_n(1+2d)$. Consequently, \citet[Eq.~(18.5.1)]{dlmf} gives that 
		$$
		\norm{\cT_n(X)}\leq \cT_n(1+2d)
		=\cosh\!\bigl(n\operatorname{arcosh}(1+2d)\bigr)\leq\ee^{2n^2d},
		$$
		where we used $\cosh t\leq\ee^{t^2/2}$ and $\operatorname{arcosh}(1+2d)\leq2\sqrt d$. Hence we complete the proof of \eqref{eq:T-growth}.
		
		For \eqref{eq:Q-diff}, applying the same identity \eqref{eq:ordered-diff} to both $I-2B$ and $I-2D$ gives 
		$$
		\begin{aligned}
			\Big\|\mathcal Q_{n,\phi}(B)&-\mathcal Q_{n,\phi}(D)\Big\|
			=\Big\|\cT_n(I-2B)-\cT_n(I-2D)\Big\|\\
			& =\Big\|-2(B-D)U_{n-1}(I-2D) -4\sum_{k=1}^{n-1} \cT_k(I-2B)(B-D)U_{n-1-k}(I-2D)\Big\|\\
			& \overset{(i)}{\leq} 2\norm{B-D}\ee^{2n^2d} \Big(n+2\sum_{k=1}^{n-1}(n-k)\Big) = 2n^2\ee^{2n^2d}\norm{B-D}, 
		\end{aligned}
		$$
		where (i) is by the following bound on the $U_k(I-2D)$ terms: we iterate $U_n(z)-U_{n-2}(z)=2\cT_n(z)$, see \citet[Eq.~(18.9.9)]{dlmf}, with $U_0=\cT_0$ and $U_1=2\cT_1$, and then applying \eqref{eq:T-growth} gives
		$$
		\norm{U_n(I-2D)} = \bigg\|\sum_{k=0}^{n}\cT_{|n-2k|}(I-2D)\bigg\|\leq\sum_{k=0}^{n}\ee^{2(n-2k)^2d} \leq(n+1)\ee^{2n^2d}.
		$$
		Finally, for $\phi\in[\pi/2,\pi)$, the definition of $\mathcal Q_{n,\phi}$ and \eqref{eq:T-growth} give
		$$
		\left\|\frac{\mathcal Q_{n,\phi}(B)}{1-\cos\phi}\right\| = \frac{\norm{\cT_n(I-2B)-\cos\phi\,I}}{1-\cos\phi} \leq \frac{\ee^{2n^2d}-\cos\phi}{1-\cos\phi} \leq \ee^{2n^2d}\leq \ee^{4n^2d},
		$$
		where the second inequality uses $\cos\phi\leq0$. This completes the proof of \eqref{eq:frozen-good-new}.
	\end{proof}
	
	\begin{lemma}\label{lem:varying-product}
		Let $\Acal_{m,\pi/2} = \big(a_1,\ldots,a_{m}\big)$ for some $m\in\{2^\ell:\ell=1,2,\cdots\}$. Fix $1\leq k<m$ and $d,W\geq0$, suppose matrices $Z_i,B_i$ satisfy $\norm{Z_i-B_i}\leq d$, $0\preceq B_i\preceq I/2$, $1\leq i\leq k$, then 
		\begin{equation} 
			\Big\|(I-a_kZ_k)\cdots(I-a_1Z_1)\Big\| \leq \exp\!\bigg(\frac43m^2d + \frac{m^2}{2}(\log_2m-1)W\bigg) \nonumber
		\end{equation}
		holds, as long as we have $\sum_{i=2}^{k}\norm{Z_i-Z_{i-1}}\leq W$. 
	\end{lemma}
	\begin{proof}
		By Lemma~\ref{lem:preceding-structure}(i), the sub-block $(a_1,\cdots,a_k)$ can be decomposed into recursive sub-blocks with phases in $[\pi/2,\pi)$. For any recursive sub-block $(a_s,\ldots,a_t)=\Acal_{m',\phi'}$, set
		$$
		W_{s:t}=\sum_{i=s+1}^{t}\norm{Z_i-Z_{i-1}}\quad\mbox{and}\quad
		\widetilde Q_{s:t}=(1-\cos\phi')(I-a_tZ_t)\cdots(I-a_sZ_s).
		$$
		Set $e_{s:t}:=\ee^{-4(m')^2d}\norm{\widetilde Q_{s:t}-\mathcal Q_{m',\phi'}(Z_t)}$, for each sub-block in this decomposition, \eqref{eq:frozen-good-new} gives
		\begin{equation}\label{eq:varying-block-bound}
			\begin{aligned}
				\norm{(I-a_tZ_t)&\cdots(I-a_sZ_s)}
				=\frac{\norm{\widetilde Q_{s:t}}}{1-\cos\phi'}\\
				&\leq \frac{\norm{\mathcal Q_{m',\phi'}(Z_t)}}{{1-\cos\phi'}} + \norm{\widetilde Q_{s:t}-\mathcal Q_{m',\phi'}(Z_t)}\\
				&\leq \ee^{4(m')^2d}(1+e_{s:t}).
			\end{aligned}
		\end{equation}
		Then it remains to control $e_{s:t}$ by induction that 
		\begin{equation}
			\label{eq:complete-product-error} 
			e_{s:t} \leq \ee^{2(m')^2 \log_2(m')\,W_{s:t}}-1.
		\end{equation} First, $e_{s:t}=0$ for $m'=1$,and hence \eqref{eq:complete-product-error} holds for $\Acal_{m',\phi'}, \forall \phi'\in(0,\pi)$.  For $m'\geq2$, assume \eqref{eq:complete-product-error} holds for all  $\Acal_{m'',\phi''}$, with $m''= m'/2$ and $\forall \phi''\in(0,\pi)$. The recursive definition \eqref{eq:A-def} and the identity $\cT_{m'}(z)=2\cT_{m''}(z)^2-1$ in \eqref{eq:cheb-standard} give
		$$
		\widetilde Q_{s:t} = 2\widetilde Q_{t-m''+1:t}\widetilde Q_{s:t-m''},\qquad
		\mathcal Q_{m',\phi'}(Z_t) =2\mathcal Q_{m'',\phi'/2}(Z_t)\mathcal Q_{m'',\pi-\phi'/2}(Z_t).
		$$
		To bound the difference of these products, we use \eqref{eq:T-growth} to obtain 
		$$
		\norm{\mathcal Q_{m'',\psi}(Z)}
		\leq\norm{\cT_{m''}(I-2Z)}+1
		\leq2\ee^{(m')^2d/2},
		\qquad Z\in\{Z_{t-m''},Z_t\},\quad \psi\in(0,\pi).
		$$
		Expanding the difference and applying \eqref{eq:Q-diff}, together with $d\geq0$, yields
		$$
		\begin{aligned}
			e_{s:t}
			&=2\ee^{-4(m')^2d}\big\|\widetilde Q_{t-m''+1:t}\widetilde Q_{s:t-m''}
			-\mathcal Q_{m'',\phi'/2}(Z_t)\mathcal Q_{m'',\pi-\phi'/2}(Z_t)\big\|\\
			&\leq4e_{s:t-m''}+4e_{t-m''+1:t}+2e_{s:t-m''}e_{t-m''+1:t}+2(m')^2\norm{Z_{t-m''}-Z_t}.
		\end{aligned}
		$$
		By the fact that  $1+4x+4y+2xy\leq(1+x)^4(1+y)^4$ for $x,y\geq0$ and $1+u\leq\ee^u$,  we obtain
		$$
		\begin{aligned}
			1+e_{s:t}
			&\leq(1+e_{s:t-m''})^4(1+e_{t-m''+1:t})^4\ee^{2(m')^2\norm{Z_{t-m''}-Z_t}}\\
			&\leq\exp\!\left(2(m')^2\big[(\log_2m'-1)(W_{s:t-m''}+W_{t-m''+1:t})+\norm{Z_{t-m''}-Z_t}\big]\right)\\
			&\overset{(i)}{\leq}\ee^{2(m')^2\log_2m'\,W_{s:t}},
		\end{aligned}
		$$
		where (i) uses $W_{s:t-m''}+W_{t-m''+1:t}\leq W_{s:t}$ and $\norm{Z_{t-m''}-Z_t}\leq W_{s:t}$. Hence \eqref{eq:complete-product-error} is true.
		
		Finally, the decomposition of $(a_1,\ldots,a_k)$ in Lemma~\ref{lem:preceding-structure}(i) contains at most one block of each length $m/2,m/4,\ldots,1$. For each such block, $(m')^2\log_2m'\leq {m^2}(\log_2m-1)/4,$
		and the $W_{s:t}$'s sum up to at most $W$. Multiplying \eqref{eq:varying-block-bound} over these sub-blocks and substituting \eqref{eq:complete-product-error} now gives
		$$
		\begin{aligned}
			\norm{(I-a_kZ_k)\cdots(I-a_1Z_1)}
			&\leq\exp\!\left(4d\sum_{\ell=1}^{\log_2m}\left(\frac{m}{2^\ell}\right)^2+\frac{m^2}{2}(\log_2m-1)W\right)\\
			&\leq\exp\!\left(\frac43m^2d+\frac{m^2}{2}(\log_2m-1)W\right).
		\end{aligned}
		$$
		Hence we complete the proof of this lemma.
	\end{proof}
	
	With the above preparation, we are ready to establish the proof of inequality \eqref{eq:positive-prod-res}, which will be presented in the next subsection. 
	
	\subsection{Proof of \eqref{eq:positive-prod-res}}
	\label{sec:positive-prod-res}
	The proof of this inequality will follow the notations of Section \ref{sec:positive-curvature-control}. Now, fix an index $i<m$. To bound $\Psi_{0,j}, j\leq i$, applying Lemma~\ref{lem:varying-product} with $Z_q={\Pi_+\overline H_q\Pi_+}/{2L}, 0\leq q<j$,  and $d=W={V_i}/{2L}$ 
	gives
	$$
	\norm{\Psi_{0,j}} = \norm{(I-a_{j-1}Z_{j-1})\cdots(I-a_0Z_0)} \leq\exp\!\bigg( \frac{m^2V_i}{2L} \left( \frac43 + \frac{\log_2m-1}{2} \right) \bigg) = \ee^{c_mV_i}.
	$$
	Here, the first equality uses $a_qZ_q=\eta_q\Pi_+\overline H_q\Pi_+$. In particular, the choice of $d=W=\frac{V_i}{2L}$ is valid because $0\preceq {\widehat H_+}/{2L}\preceq I/2$,
	and by \eqref{eq:epoch-hessian-variation}, we have
	$$
	\begin{aligned}
		\left\|Z_q - {\widehat H_+}/{2L}\right\|
		&=\frac{\norm{\Pi_+(\overline H_q-\widehat H)\Pi_+}}{2L} \leq \frac{\norm{\overline H_q-\widehat H}}{2L} \leq\frac{V_i}{2L},
		\qquad 0\leq q<j,\\[4pt]
		\sum_{q=1}^{j-1}\norm{Z_q-Z_{q-1}}
		&\leq\frac1{2L}
		\sum_{q=1}^{j-1}\norm{\overline H_q-\overline H_{q-1}}
		\leq\frac{V_i}{2L}.
	\end{aligned}
	$$
	This justifies our use of Lemma~\ref{lem:varying-product} here. Next, to bound $F_{\ell,j}$ for $0\leq\ell<j\leq i$, we represent its action on a unit vector as a block of a matrix product. For any unit vector $u\in\mathbb R^d$ and $\rho>0$, denote $\mathbf 1_{\{q\geq\ell\}}$ the indicator that takes value 1 when $q\geq\ell$ and $0$ otherwise, and define
	$$
	Z_q^{(\rho)}
	:=\frac1{2L}
	\begin{pmatrix}
		\Pi_+\overline H_q\Pi_+ & \rho\mathbf 1_{\{q\geq\ell\}}u\\
		0 & 0
	\end{pmatrix},
	\qquad 0\leq q<j.
	$$
	Then these $Z_q^(\rho)$'s are not symmetric and not commutative. Multiplying these blocks gives 
	$$
	\big(I-a_{j-1}Z_{j-1}^{(\rho)}\big)\cdots\big(I-a_0Z_0^{(\rho)}\big)
	=
	\begin{pmatrix}
		\Psi_{0,j} & -\rho F_{\ell,j}u\\
		0 & 1
	\end{pmatrix}.
	$$
	Because each $Z_q^{(\rho)}$ is within  $\frac{V_i+\rho}{2L}$ distance of symmetric matrix $0\preceq \frac{1}{2L}\mathrm{Diag}(\widehat H_+, 0)\preceq\frac{I}{2}$, and $\sum_{q=1}^{j-1}\norm{Z_q^{(\rho)}-Z_{q-1}^{(\rho)}}\leq \frac{V_i+\rho}{2L}$ since the indicator $\mathbf 1_{\{q\geq\ell\}}$ changes at most once, we get 
	$$ 
	\left\|(I-a_{j-1}Z_{j-1}^{(\rho)})\cdots(I-a_0Z_0^{(\rho)})\right\| \leq \exp\!\bigg(
	\frac{m^2(V_i+\rho)}{2L}\left(\frac43 +\frac{\log_2m-1}{2}\right)\bigg]
	= \ee^{c_m(V_i+\rho)}, 
	$$
	by applying Lemma~\ref{lem:varying-product}. Therefore, we have  
	$$
	\rho\norm{F_{\ell,j}u} \leq 
	\left\|\begin{pmatrix}
		\Psi_{0,j} & -\rho F_{\ell,j}u\\
		0 & 1
	\end{pmatrix}\right\|
	\leq \ee^{c_m(V_i+\rho)}.
	$$
	Finally, set $\rho=1/c_m$ and take  supremum over $\{u:\|u\|=1\}$ yields
	$\norm{F_{\ell,j}}\leq c_m\ee^{1+c_mV_i}.$
	Hence we complete the proof of  \eqref{eq:positive-prod-res}.
	
	\section{Proof of Lemma~\ref{lem:convex-local-safety}}
	\label{sec:convex-technical}
	\begin{proof}
		Fix an epoch $t$ and an iterate $i$ with $\eta_{t,i}>1/L$. Similar to  Section~\ref{sec:positive-curvature-control}, we omit the epoch index $t$ and write $x_j=x_{t,j}$, $g_j=\nabla f(x_j)$, $a_j=a_{t,j}$, and $G=G_t$. Let $H_i:=\nabla^2f(x_i)$, and use $\overline H_j$ for the averaged Hessian in \eqref{eq:epoch-gradient-dynamics}. As clipping ends an epoch,  $\eta_j=\frac{a_j}{2L}$ for $j<i$.

		First of all, let us control the cubic term of this lemma. As Lemma~\ref{lem:nonconvex-path} still holds, we have $\eta_i\norm{g_i}\leq\frac{\theta p_mG}{4L}$. Together with $\eta_i\leq p_m/2L$ and \eqref{eq:nonconvex-cert}, we obtain
		\begin{equation}
			\label{eq:cubic-convex}
			\frac{M\eta_i^2}{3}\norm{g_i}^3 = \frac{M}{3}\cdot \eta_i\norm{g_i} \cdot \eta_i\cdot \norm{g_i}^2 \leq \frac{\theta MGp_m^2}{24L^2}\norm{g_i}^2
			\leq\frac{\theta}{72h_m^2}\norm{g_i}^2.
		\end{equation}
		Next, we bound the quadratic term of this lemma. By Lemma~\ref{lem:preceding-structure}(ii), let $Q=[a_s,\ldots,a_i] = \Acal_{m',\phi'}$ be the largest recursive sub-block ending at $a_i$, with $\phi\in[\pi/2,\pi)$. Since $a_i\geq2L\eta_i>2$, the same lemma implies $m'\geq2$.  Since $a_i$ is the last coefficient of $Q$, Lemma~\ref{lem:dissipative} gives $(1+a_i z/6)F_Q(z)^2\leq1$ for $z\in[0,1/2]$. Applying this   inequality to $0\preceq \frac{1}{2L}H_i\preceq \frac{1}{2}I$ gives 
		$$
		\norm{g_s}^2 \geq g_s^\top\left(I+\frac{a_i}{6}\cdot\frac{H_i}{2L}\right) F_Q\!\left(\frac{H_i}{2L}\right)^2g_s = \norm{w_i}^2+\frac16\bigg\|\sqrt{\frac{a_iH_i}{2L}}\,w_i\bigg\|^2
		$$
		where $w_i:=F_Q\big(\frac{1}{2L}H_i\big)g_s$.  Consequently,  $\eta_i\leq a_i/2L$ and triangle inequality together gives   
		\begin{equation}\label{eq:convex-curvature-preparation}
			\begin{aligned}
				\sqrt{\eta_i g_i^\top H_i g_i}\leq \bigg\|\sqrt{\frac{a_iH_i}{2L}}\,g_i\bigg\|
				&\leq\bigg\|\sqrt{\frac{a_iH_i}{2L}}\,w_i\bigg\|
				+\bigg\|\sqrt{\frac{a_iH_i}{2L}}\,(g_i-w_i)\bigg\|\\
				&\leq\sqrt{6\bigl(\norm{g_s}^2-\norm{w_i}^2\bigr)}
				+\bigg\|\sqrt{\frac{a_iH_i}{2L}}\,(g_i-w_i)\bigg\|.
			\end{aligned}
		\end{equation}\vspace{0.2cm}
		Next, we bound the above terms one by one. 
		
		\noindent\textbf{Bounding the gradients.}\,\,
		By convexity, setting $\Pi_+=I$ and $\Pi_-=0$ in \eqref{eq:positive-forced-stability} gives
		\begin{equation}\label{eq:convex-envelope}
			\norm{g_j}\leq\ee^{c_mV_i}G
			\leq\ee^\nu G,\qquad 0\leq j\leq i.
		\end{equation}
		Here the second inequality is by \eqref{eq:positive-gradient-exponent}, with $\nu=\frac{\theta m^2}{24p_m}\left(\frac{1}{2h_m}+\frac{1}{3h_m^2}\right)$ defined in Lemma~\ref{lem:positive-gradient-bound}.\vspace{0.2cm}
		
		\noindent\textbf{Bounding $\norm{g_i-w_i}$ and $\norm{\sqrt{\frac{a_iH_i}{2L}}\,(g_i-w_i)}$.}\,\,
		First, rewrite gradient dynamics \eqref{eq:epoch-gradient-dynamics} as 
		$$g_{j+1}=(I-\eta_jH_i)g_j+\eta_j(H_i-\overline H_j)g_j.$$
		Unrolling from $s$ to $i$ and using $\eta_j=a_j/2L$ gives
		\begin{equation}\label{eq:convex-error-expansion}
			g_i-w_i=\sum_{\ell=1}^{m'-1}\frac{a_{s+\ell-1}}{2L}
			P_{Q_{\ell+1:m'-1}}\!\left(\frac{H_i}{2L}\right)
			(H_i-\overline H_{s+\ell-1})g_{s+\ell-1}.
		\end{equation}
		By Lemma~\ref{lem:tail-preparation}, each polynomial term in this summation has norm at most 1. To bound the Hessian differences, we can use Hessian $M$-Lipschitz continuity and Lemma~\ref{lem:nonconvex-path} to obtain
		\begin{equation}\label{eq:convex-hessian-distance}
			\max_{0\leq j<i}\norm{H_i-\overline H_j}
			\leq M\sum_{j=0}^{i-1}\eta_j\norm{g_j}
			\leq\frac{\theta Mp_mG}{4L}.
		\end{equation}
		Substituting into \eqref{eq:convex-error-expansion} and applying  Lemma~\ref{lem:nonconvex-path} and \eqref{eq:nonconvex-cert}, we obtain
		\begin{align*}
			\norm{g_i-w_i} \leq\sum_{j=s}^{i-1}\eta_j\norm{H_i-\overline H_j}\norm{g_j} \leq M\bigg(\sum_{j=0}^{i-1}\eta_j\norm{g_j}\bigg)^2 \leq\frac{\theta^2Mp_m^2G^2}{16L^2}
			\leq\frac{\theta^2}{48h_m^2}G.
		\end{align*} 
		Next, applying \eqref{eq:convex-envelope} and \eqref{eq:convex-hessian-distance} to \eqref{eq:convex-error-expansion} gives
		\begin{align*}
			\bigg\|\sqrt{\frac{a_i}{2L}H_i} (g_i-w_i)\bigg\| &\leq \ee^\nu\frac{\theta Mp_mG^2}{8L^2}\!
			\sum_{\ell=1}^{m'-1}\!\!a_{s+\ell-1}
			\left\|\sqrt{\frac{a_i}{2L}H_i} P_{Q_{\ell+1:m'-1}}\!\left(\frac{H_i}{2L}\right)\right\|\\
			&\leq \frac{\theta\pi\ee^\nu m^2(h_m\!\!-\!1)}{24\sqrt{2}p_mh_m^2} G 
		\end{align*}
		where the last inequality is by applying Lemma~\ref{lem:weighted-response} to $Q$ and $H_i/2L$, with $m'\leq m$ and \eqref{eq:nonconvex-cert}.\vspace{0.2cm}

		\noindent\textbf{Finalizing the proof.}\,\,
		Combining $\norm{g_i}>\theta G$ and the above bound on $\norm{g_i-w_i}$, we have 
		$$\norm{w_i}\geq\norm{g_i}-\norm{g_i-w_i} > \big(\theta-\frac{\theta^2}{48h_m^2}\big)G.$$ 
		Substituting this inequality, \eqref{eq:convex-envelope}, and the bound on $\norm{\sqrt{\frac{a_iH_i}{2L}}\,(g_i-w_i)}$ into \eqref{eq:convex-curvature-preparation} gives 
		\begin{equation}
			\eta_i g_i^\top H_i g_i \leq\left[\sqrt{6\left\{\ee^{2\nu}
				-\left(\theta-\frac{\theta^2}{48h_m^2}\right)^2\right\}}
			+\frac{\theta\pi\ee^\nu m^2(h_m-1)}{24\sqrt2\,p_mh_m^2}\right]^2\!\!G^2 \leq\left(1-\frac{\theta}{72h_m^2}\right)\norm{g_i}^2. \nonumber
		\end{equation}
		Here the last inequality uses $\norm{g_i}>\theta G$, $\theta=19/20$, $h_m\geq3$, $m^2/p_m\leq\pi^2/16$, and the definition of $\nu$. Combining this bound with \eqref{eq:cubic-convex} gives
		$$
		\eta_i g_i^\top H_i g_i+\frac{M\eta_i^2}{3}\norm{g_i}^3
		\leq\norm{g_i}^2,
		$$
		which completes the proof of Lemma~\ref{lem:convex-local-safety}.
	\end{proof}

	\bibliography{ref}
	
\end{document}